\documentclass{article}

\usepackage{amsthm}
\usepackage{amsmath}
\usepackage{amssymb}
\usepackage{mathrsfs}
\usepackage{dsfont}
\usepackage[all]{xy}
\usepackage{tikz-cd}
\usepackage{adjustbox}
\usepackage{enumitem}
\usepackage{newclude}
\usepackage{hyperref}

\theoremstyle{definition}
\newtheorem{Definition}{Definition}[section]

\theoremstyle{plain}
\newtheorem{Theorem}[Definition]{Theorem}

\theoremstyle{plain}
\newtheorem{Proposition}[Definition]{Proposition}

\theoremstyle{plain}
\newtheorem{Lemma}[Definition]{Lemma}

\theoremstyle{plain}
\newtheorem{Corollary}[Definition]{Corollary}

\theoremstyle{definition}
\newtheorem{Example}[Definition]{Example}

\theoremstyle{remark}
\newtheorem{Remark}[Definition]{Remark}

\theoremstyle{plain}
\newcommand{\thistheoremname}{}
\newtheorem*{genericthm*}{\thistheoremname}
\newenvironment{namedthm*}[1]
  {\renewcommand{\thistheoremname}{#1}%
   \begin{genericthm*}}
  {\end{genericthm*}}

\author{Thibault D. Décoppet}
\title{The 2-Deligne Tensor Product}

\newenvironment{tz}[1][]{\begin{aligned}\begin{tikzpicture}[#1]}{\end{tikzpicture}\end{aligned}}
\tikzset{morphlabel/.style={draw=black, thin, rectangle, minimum width=7pt, fill=white, font=\scriptsize}}

\newcommand*\arrowoffset{0.4pt}
\makeatletter
\pgfarrowsdeclare{new double arrowhead}{new double arrowhead}
{
  \pgfutil@tempdima=-0.84pt%
  \advance\pgfutil@tempdima by-1.3\pgflinewidth%
  \pgfutil@tempdimb=0.21pt%
  \advance\pgfutil@tempdimb by.625\pgflinewidth%
  \pgfarrowsleftextend{+\pgfutil@tempdima}
  \advance\pgfutil@tempdimb by \arrowoffset % Jamie addition
  \pgfarrowsrightextend{+\pgfutil@tempdimb}
}
{
  \pgfsetstrokecolor{black}
  \pgfsetlinewidth{0.6pt}
  \pgfutil@tempdima=0.28pt%
  \advance\pgfutil@tempdima by.3\pgflinewidth%
  \pgfsetlinewidth{0.8\pgflinewidth}
  \pgfsetdash{}{+0pt}
  \pgfsetroundcap
  \pgfsetroundjoin
  \pgfpathmoveto{\pgfpoint{-3\pgfutil@tempdima+\arrowoffset}{4\pgfutil@tempdima}}
  \pgfpathcurveto
                {\pgfpoint{-2.75\pgfutil@tempdima+\arrowoffset}{2.5\pgfutil@tempdima}}
                {\pgfpoint{\arrowoffset}{0.25\pgfutil@tempdima}}
                {\pgfpoint{0.75\pgfutil@tempdima+\arrowoffset}{0pt}}
  \pgfpathcurveto
                {\pgfpoint{\arrowoffset}{-0.25\pgfutil@tempdima}}
                {\pgfpoint{-2.75\pgfutil@tempdima+\arrowoffset}{-2.5\pgfutil@tempdima}}
                {\pgfpoint{-3\pgfutil@tempdima+\arrowoffset}{-4\pgfutil@tempdima}}
  \pgfusepathqstroke
}
\pgfarrowsdeclare{new single arrowhead}{new single arrowhead}
{
  \pgfutil@tempdima=-0.84pt%
  \advance\pgfutil@tempdima by-1.3\pgflinewidth%
  \pgfutil@tempdimb=0.21pt%
  \advance\pgfutil@tempdimb by.625\pgflinewidth%
  \pgfarrowsleftextend{+\pgfutil@tempdima}
  \advance\pgfutil@tempdimb by \arrowoffset % Jamie addition
  \pgfarrowsrightextend{+\pgfutil@tempdimb}
}
{
  \pgfsetstrokecolor{black}
  \pgfsetlinewidth{0.6pt}
  \pgfutil@tempdima=0.28pt%
  \advance\pgfutil@tempdima by.3\pgflinewidth%
  \pgfsetlinewidth{0.8\pgflinewidth}
  \pgfsetdash{}{+0pt}
  \pgfsetroundcap
  \pgfsetroundjoin
  \pgfpathmoveto{\pgfpoint{-3\pgfutil@tempdima}{4\pgfutil@tempdima}}
  \pgfpathcurveto
                {\pgfpoint{-2.75\pgfutil@tempdima}{2.5\pgfutil@tempdima}}
                {\pgfpoint{0pt}{0.25\pgfutil@tempdima}}
                {\pgfpoint{0.75\pgfutil@tempdima}{0pt}}
  \pgfpathcurveto
                {\pgfpoint{0pt}{-0.25\pgfutil@tempdima}}
                {\pgfpoint{-2.75\pgfutil@tempdima}{-2.5\pgfutil@tempdima}}
                {\pgfpoint{-3\pgfutil@tempdima}{-4\pgfutil@tempdima}}
  \pgfusepathqstroke
}
\makeatother

\tikzset{morphlabel/.style={draw=black, thin, rectangle, minimum width=7pt, fill=white, font=\scriptsize}}

\tikzset{double arrow scope/.style={every path/.style={double, -new double arrowhead}}}

\begin{document}

\bibliographystyle{alpha}

    \maketitle
    \hspace{1cm}
    \begin{abstract}
        We prove that the 2-Deligne tensor product of two compact semisimple 2-categories exists. Further, under suitable hypotheses, we explain how to describe the $Hom$-categories, connected components, and simple objects of a 2-Deligne tensor product. Finally, we prove that the 2-Deligne tensor product of two compact semisimple tensor 2-categories is a compact semisimple tensor 2-category.
    \end{abstract}
    
\tableofcontents
    
\section*{Introduction}
\addcontentsline{toc}{section}{Introduction}

In \cite{D2}, we began our investigation of the algebraic properties of finite semisimple 2-categories and multifusion 2-categories over an algebraically closed field of characteristic zero (see also \cite{DR}, and \cite{JF}). Then, in \cite{D5}, working over an arbitrary field, we defined compact semisimple 2-categories, which generalize finite semisimple 2-categories. Namely, whereas finite semisimple 2-categories have finitely many equivalence classes of simple objects, compact semisimple 2-categories are merely required to have finitely many connected components. Let us note that over an algebraically closed field or a real closed field, every compact semisimple 2-category is finite (see \cite{D5}). Over arbitrary fields however, it becomes necessary to introduce the notion of a compact semisimple 2-category because finite semisimple 2-categories do not exist in general. Thus, compact semisimple 2-categories are the appropriate categorification of the notion of a finite semisimple category over an arbitrary field. Similarly, compact semisimple tensor 2-categories categorify finite semisimple tensor categories. Presently, we continue our investigations of the properties of compact semisimple 2-categories and compact semisimple tensor 2-categories by describing a 2-categorical analogue of the so-called Deligne tensor product of finite categories over perfect fields (see \cite{Del}, \cite{EGNO}, and \cite{DSPS13}).

We start in section \ref{sec:1Deligne} by recalling the definition of the (1-)Deligne tensor product. Namely, given  $\mathcal{C}$ and $\mathcal{D}$ two finite linear categories (see \cite{EGNO} for the definition), their (1-)Deligne tensor product $\boxtimes:\mathcal{C}\times\mathcal{D}\rightarrow\mathcal{C}\boxtimes\mathcal{D}$ is the 2-universal bilinear functor to a finite linear category, which is right-exact in both variables. If $\mathcal{C}$ and $\mathcal{D}$ are finite semisimple categories, we can use a completed tensor product to compute their Deligne tensor product. Namely, for arbitrary linear categories $\mathcal{A}$ and $\mathcal{B}$, one can define their tensor product $\mathcal{A}\otimes\mathcal{B}$ by taking the usual tensor product of their $Hom$-spaces. Additionally, we can Cauchy complete $\mathcal{A}\otimes\mathcal{B}$, i.e. formally add direct sums and splittings of idempotents, and get a linear category $Cau(\mathcal{A}\otimes \mathcal{B})$. It turns out that if we take as input for this completed tensor product two finite semisimple categories $\mathcal{C}$ and $\mathcal{D}$, then $Cau(\mathcal{C}\otimes \mathcal{D})$ is finite semisimple, provided the base field is perfect. This implies that the 2-universal bilinear functor $\mathcal{C}\times\mathcal{D}\rightarrow Cau(\mathcal{C}\otimes \mathcal{D})$ coincides with the Deligne tensor product $\boxtimes:\mathcal{C}\times\mathcal{D}\rightarrow\mathcal{C}\boxtimes\mathcal{D}$. It is this second approach we will use in order to construct the desired 2-Deligne tensor products. For later use, we also recall how to compute certain $Hom$-spaces in a Deligne tensor product. More precisely, if $C_1$, $C_2$ are in $\mathcal{C}$ and $D_1$, $D_2$ are in $\mathcal{D}$, we have $$Hom_{\mathcal{C}\boxtimes\mathcal{D}}(C_1\boxtimes D_1, C_2\boxtimes D_2)\cong Hom_{\mathcal{C}}(C_1,C_2)\otimes Hom_{\mathcal{D}}(D_1,D_2).$$ We also review the key fact that, over a perfect field, the Deligne tensor product of two finite tensor categories is a finite tensor category.

Section \ref{sec:recollections} begins by recalling the 3-universal property of the Cauchy completion of a linear 2-category from \cite{D1}. (The notion of Cauchy completion of a higher category was introduced in \cite{GJF}.) We move on to review the definitions of a semisimple 2-category and of a finite semisimple 2-category over an arbitrary field (see \cite{DR} for the definition over an algebraically closed field of characteristic zero). Moreover, we recall the definition of a compact semisimple 2-category from \cite{D5}, which is the appropriate categorification of the notion of a finite semisimple category. We review the crucial fact that every compact semisimple 2-category is equivalent to the 2-category of separable module categories over a finite semisimple tensor category.

Then, in section \ref{sec:2Deligne}, given two linear 2-categories $\mathfrak{A}$ and $\mathfrak{B}$, we go on to define their completed tensor product $\widehat{\otimes}:\mathfrak{A}\times\mathfrak{B}\rightarrow \mathfrak{A}\widehat{\otimes}\mathfrak{B}$, as the 3-universal bilinear 2-functor from $\mathfrak{A}\times\mathfrak{B}$ to a Cauchy complete 2-category, and we show that it always exists. Let us now assume that our base field is perfect. Given $\mathfrak{C}$ and $\mathfrak{D}$ two compact semisimple 2-categories, we define their 2-Deligne tensor product $\boxdot:\mathfrak{C}\times\mathfrak{D}\rightarrow\mathfrak{C}\boxdot\mathfrak{D}$, if it exists, as the 3-universal bilinear 2-functor from $\mathfrak{C}\times \mathfrak{D}$ to a compact semisimple 2-category. By analogy with the decategorified situation, the 2-Deligne tensor product can be computed using the completed tensor product $\widehat{\otimes}:\mathfrak{C}\times\mathfrak{D}\rightarrow \mathfrak{C}\widehat{\otimes}\mathfrak{D}$, provided that $\mathfrak{C}\widehat{\otimes}\mathfrak{D}$ is a compact semisimple 2-category. This is precisely how the following result is proven.

\begin{namedthm*}{Theorem \ref{thm:2Deligne}}
Over a perfect field, the 2-Deligne tensor product $\boxdot:\mathfrak{C}\times\mathfrak{D}\rightarrow\mathfrak{C}\boxdot\mathfrak{D}$ of two compact semisimple 2-categories $\mathfrak{C}$ and $\mathfrak{D}$ exists.
\end{namedthm*}

\noindent Let us mention that, over an algebraically closed field of characteristic zero, the idea of constructing the $n$-Deligne tensor product of two separable finite semisimple $n$-categories using a suitably completed tensor products is also featured in \cite{JF}. However, it should be noted that our setup is more general than theirs in the case $n=2$, given that we work over an arbitrary perfect field. Moreover, it is essential to understand the 2-Deligne tensor product of two multifusion 2-categories. Yet, in \cite{JF}, they only examine the 2-Deligne tensor product of two separable multifusion 2-categories, and it remains an open problem to prove that, over an algebraically closed field of characteristic zero, every multifusion 2-category is separable. We believe that our approach to the 2-Deligne tensor product of compact semisimple 2-categories will play a role in the proof of this result.

In section \ref{sec:properties}, we examine the properties of the 2-Deligne tensor product of two compact semisimple 2-categories. Some of these features resemble that of the (1-)Deligne tensor product. In particular, given objects $C_1,C_2$ in $\mathfrak{C}$, and $D_1,D_2$ in $\mathfrak{D}$, we show that there is an equivalence $$Hom_{\mathfrak{C}\boxdot\mathfrak{D}}(C_1\boxdot D_1, C_2\boxdot D_2)\simeq Hom_{\mathfrak{C}}(C_1, C_2)\boxtimes Hom_{\mathfrak{D}}(D_1, D_2).$$ We wish to emphasize that, even when working over an algebraically closed field of characteristic zero, the behaviour of the 2-Deligne tensor product is closer to that of the (1-)Deligne tensor product over an arbitrary perfect field. For instance, it may happen that the equivalence classes of simple objects of $\mathfrak{C}\boxdot\mathfrak{D}$ are not parametrised by pairs consisting of an equivalence class of simple objects in $\mathfrak{C}$, and one in $\mathfrak{D}$. We give an example of this behaviour. On the other hand, over an algebraically closed field, we can define the Frobenius-Perron dimension of a connected finite semisimple 2-category $\mathfrak{C}$ as the Frobenius-Perron dimension of the fusion category $End_{\mathfrak{C}}(C)$ for any simple object $C$ of $\mathfrak{C}$. Note that this is an algebraic integer. Using this definition, we can prove the following result:

\begin{namedthm*}{Theorem \ref{thm:2Deligneessentiallysurjective}}
Let $\mathfrak{C}$ and $\mathfrak{D}$ be connected finite semisimple 2-categories over an algebraically closed field. If $\mathfrak{C}$ and $\mathfrak{D}$ have coprime Frobenius-Perron dimension, then, for any simple object $A$ of $\mathfrak{C}\boxdot\mathfrak{D}$, there exists simple objects $C$ in $\mathfrak{C}$ and $D$ in $\mathfrak{D}$ unique up to equivalence such that $C\boxdot D \simeq A$.
\end{namedthm*}

Finally, in section \ref{sec:2DeligneMonoidal}, we show that the completed tensor $\mathfrak{A}\widehat{\otimes}\mathfrak{B}$ of two monoidal linear categories $\mathfrak{A}$ and $\mathfrak{B}$ inherits a monoidal structure. In particular, over a perfect field, the 2-Deligne tensor product can be used to create new compact semisimple tensor 2-categories, i.e. compact semisimple 2-categories equipped with a rigid monoidal structure, out of the ones we already know.

\begin{namedthm*}{Theorem \ref{thm:2delignetensor}}
Given $\mathfrak{C}$ and $\mathfrak{D}$ two compact semisimple tensor 2-categories over a perfect field. Their 2-Deligne tensor product $\mathfrak{C}\boxdot\mathfrak{D}$ is a compact semisimple tensor 2-category. Further, the 2-functor $\boxdot:\mathfrak{C}\times\mathfrak{D}\rightarrow \mathfrak{C}\boxdot\mathfrak{D}$ is monoidal.
\end{namedthm*}

\noindent We end by giving some examples of the fusion rule of the 2-Deligne tensor product of two compact semisimple tensor 2-categories over various algebraically closed fields.

Generalising the notion of the (1-)Deligne tensor product, there is a notion of relative, or balanced, (1-)Deligne tensor product of two module categories over a fixed finite tensor category. At this point, it is therefore natural to ask what can be said about the existence of the relative 2-Deligne tensor product. By mimicking the construction of the (1-)Deligne tensor product given in \cite{DSPS14} and using the main result of \cite{D4}, it is possible to define a candidate for the relative 2-Deligne tensor product of two compact semisimple module 2-categories. However, it is not clear whether the 2-category so constructed satisfies the required 3-universal property in general. (For connected compact semisimple tensor 2-categories, this is done in \cite{BJS}.) We plan on investigating this point further in the future.

\subsubsection*{Acknowledgments}

I am very much indebted to Christopher Douglas and Andr\'e Henriques regarding the content of this article. Moreover, I thank the referee for suggesting generalizing the definition of finite semisimple 2-categories over an algebraically closed field of characteristic zero to an arbitrary field, which lead me to write \cite{D5}, and greatly improved the present article.

\section{The Deligne Tensor Product}\label{sec:1Deligne}

Let us fix $\mathds{k}$ a field. We begin by recalling a few definitions from \cite{DSPS13}.

\begin{itemize}
    \item A $\mathds{k}$-linear category is called finite if it is equivalent to the category of finite dimensional modules over a finite $\mathds{k}$-algebra.
    \item A finite category is semisimple if every object splits as a finite direct sum of simple objects.
    \item A tensor category is a rigid monoidal linear category.
\end{itemize}

\noindent We are now ready to give the definition of the (1-)Deligne tensor product.

\begin{Definition}\label{def:1Deligne}
Given $\mathcal{C}$ and $\mathcal{D}$ two finite categories, their Deligne tensor product, if it exists, is the bilinear bifunctor $\boxtimes: \mathcal{C}\times\mathcal{D}\rightarrow \mathcal{C}\boxtimes\mathcal{D}$, right-exact in both variables, satisfying the following 2-universal property:
\begin{enumerate}
\item For any finite category $\mathcal{E}$ and every bilinear bifunctor right-exact in both variables $F:\mathcal{C}\times\mathcal{D}\rightarrow \mathcal{E}$, there exists a right-exact functor $F':\mathcal{C}\boxtimes\mathcal{D}\rightarrow\mathcal{E}$ and a natural isomorphism $u:F'\circ\boxtimes\Rightarrow F$.

\item Moreover, for any right-exact functors $G,H:\mathcal{C}\boxtimes\mathcal{D}\rightarrow\mathcal{E}$ and natural transformation $t:G\circ\boxtimes\Rightarrow H\circ \boxtimes$, there exists a unique natural transformation $t':G\Rightarrow H$ such that $t' \circ \boxtimes=t$.
\end{enumerate}
\end{Definition}

\begin{Remark}\label{rem:1DeligneReformulation}
We can reformulate definition \ref{def:1Deligne} as follows: Given $\mathcal{C}$ and $\mathcal{D}$ two finite categories, their Deligne tensor product, denoted by $$\boxtimes: \mathcal{C}\times\mathcal{D}\rightarrow \mathcal{C}\boxtimes\mathcal{D},$$ is the bilinear bifunctor of finite categories right-exact in both variables such that precomposition with $\boxtimes$ induces a equivalence of categories $$Hom(\mathcal{C}\boxtimes \mathcal{D},\mathcal{E})\simeq Hom_{bil}(\mathcal{C}\times \mathcal{D},\mathcal{E}),$$ from the category of right exact functors $\mathcal{C}\boxtimes \mathcal{D}\rightarrow\mathcal{E}$ to the category of bilinear bifunctor $\mathcal{C}\times \mathcal{D}\rightarrow\mathcal{E}$ right-exact in both variables. This reformulation is well-known (for instance, see \cite{DSPS14}), whence we have the following result.
\end{Remark}

\begin{Proposition}\label{prop:1Deligne}
The Deligne tensor product of two finite categories exists, is a finite category, and is unique up to a unique natural isomorphism.
\end{Proposition}

In order to motivate our construction of the 2-Deligne tensor product of two compact semisimple 2-categories, we will now explain how the Deligne tensor product of two finite semisimple categories over a perfect field can be constructed using a completed tensor product. As the latter construction is quite general, let us now consider $\mathcal{A}$ and $\mathcal{B}$ two $\mathds{k}$-linear categories. We begin by defining the category $\mathcal{A}\otimes\mathcal{B}$ whose set of objects is $Ob(\mathcal{A})\times Ob(\mathcal{B})$, and for which the space of morphisms from $(A_1,B_1)$ to $(A_2,B_2)$ is $$Hom_{\mathcal{A}\otimes\mathcal{B}}((A_1,B_1),(A_2,B_2)):=Hom_{\mathcal{A}}(A_1,A_2)\otimes_{\mathds{k}}Hom_{\mathcal{B}}(B_1,B_2).$$ The universal property of the tensor product ensures that this is a category. Further, it has the following 2-universal property.

\begin{Lemma}\label{lem:tensorproductcategories}
Given $\mathcal{A}$, and $\mathcal{B}$ be two $\mathds{k}$-linear categories. There is a 2-universal bilinear functor $\beta:\mathcal{A}\times \mathcal{B}\rightarrow \mathcal{A}\otimes\mathcal{B}$, i.e. for every $\mathds{k}$-linear category $\mathcal{E}$, precomposition with $\beta$ induces an equivalence of categories
$$Hom(\mathcal{A}\otimes \mathcal{B},\mathcal{E})\simeq Hom_{bil}(\mathcal{A}\times \mathcal{B},\mathcal{E})$$ from the category of linear functors $\mathcal{A}\otimes \mathcal{B}\rightarrow\mathcal{E}$ to the category of bilinear functors $\mathcal{A}\times \mathcal{B}\rightarrow\mathcal{E}$.
\end{Lemma}
\begin{proof}
The universal property of the tensor product proves that the canonical bilinear functor $\beta:\mathcal{A}\times \mathcal{B}\rightarrow \mathcal{A}\otimes \mathcal{B}$ is 2-universal. 
\end{proof}

The second ingredient in the construction of the completed tensor product of two $\mathds{k}$-linear categories is the notion of Cauchy completion, i.e. completion under direct sums and splittings of idempotents. By lemma A.1.8 and proposition A.1.9 of \cite{D1}, Cauchy completion has the following 2-universal property.

\begin{Lemma}\label{lem:1CauchyCompletion}
Let $\mathcal{A}$ be a $\mathds{k}$-linear category. Then, there exists a 2-universal linear functor $\kappa:\mathcal{A}\rightarrow Cau(\mathcal{A})$ to a Cauchy complete category, i.e for every Cauchy complete $\mathds{k}$-linear category $\mathcal{E}$, precomposition by $\kappa$ induces an equivalence of categories $$Hom(Cau(\mathcal{A}),\mathcal{E})\simeq Hom(\mathcal{A},\mathcal{E})$$ between the categories of linear functors.
\end{Lemma}

Putting lemma \ref{lem:tensorproductcategories} and \ref{lem:1CauchyCompletion} together, we find that for every $\mathds{k}$-linear categories $\mathcal{A}$ and $\mathcal{B}$, the bilinear bifunctor $\mathcal{A}\times \mathcal{B}\rightarrow Cau(\mathcal{A}\otimes \mathcal{B})$ is 2-universal amongst bilinear bifunctors from $\mathcal{A}\times \mathcal{B}$ to a Cauchy complete $\mathds{k}$-linear category. If we now take as input two finite semisimple categories $\mathcal{C}$ and $\mathcal{D}$, more can be said.

\begin{Proposition}
Let $\mathcal{C}$ and $\mathcal{D}$ be two finite semisimple $\mathds{k}$-linear categories, then the completed tensor product $Cau(\mathcal{C}\otimes \mathcal{D})$ is a finite semisimple category. In particular, it agrees with the Deligne tensor product $\mathcal{C}\boxtimes\mathcal{D}$.
\end{Proposition}

\begin{proof}
Let $\mathcal{E}$ be an arbitrary finite category, and $F:\mathcal{C}\times \mathcal{D}\rightarrow \mathcal{E}$ be a bilinear bifunctor. (It is necessarily right-exact in both variables as $\mathcal{C}$ and $\mathcal{D}$ are semisimple). Thanks to lemmas \ref{lem:tensorproductcategories} and \ref{lem:1CauchyCompletion}, $F$ can be factored up to natural isomorphism through $\beta$ and $\kappa$: $$\begin{tikzcd}[sep=small]
\mathcal{C}\times\mathcal{D} \arrow[rd, "\beta"'] \arrow[rrrr, "F"] &                                                                                &                                                        &  & \mathcal{E}. \\
                                                                    & \mathcal{C}\otimes\mathcal{D} \arrow[rd, "\kappa"'] \arrow[rrru, dotted] &                                                        &  &             \\
                                                                    &                                                                                & Cau(\mathcal{C}\otimes\mathcal{D}) \arrow[rruu, dotted, "F'"'] &  &            
\end{tikzcd}$$
The functor $F'$ is 2-universal by construction. Hence, if we can show that $Cau(\mathcal{C}\otimes \mathcal{D})$ is a finite category, we will get that the composite $$\mathcal{C}\times \mathcal{D}\rightarrow \mathcal{C}\otimes \mathcal{D}\rightarrow Cau(\mathcal{C}\otimes \mathcal{D})$$ is the Deligne tensor product of $\mathcal{C}$ and $\mathcal{D}$. In fact, we will now show that $Cau(\mathcal{C}\otimes \mathcal{D})$ is finite semisimple.

Pick objects $C$ in $\mathcal{C}$ and $D$ in $\mathcal{D}$ such that the corresponding fully faithful inclusions $\kappa_C:\mathrm{B}End_{\mathcal{C}}(C)\hookrightarrow \mathcal{C}$ and $\kappa_D:\mathrm{B}End_{\mathcal{D}}(D)\hookrightarrow \mathcal{D}$ are Cauchy completions. This is possible because $\mathcal{C}$ and $\mathcal{D}$ are finite semisimple categories. Firstly, observe that $$\mathrm{B}End_{\mathcal{C}}\otimes \mathrm{B}End_{\mathcal{D}} \cong \mathrm{B}(End_{\mathcal{C}}(C)\otimes End_{\mathcal{D}}(D))$$ by construction. Secondly, note that the image of $\kappa_C\otimes \kappa_D$ generates $\mathcal{C}\otimes\mathcal{D}$ under direct sums and idempotent splittings. Using lemma A.2.4 of \cite{D1}, this implies that $$Cau(\mathrm{B}(End_{\mathcal{C}}(C)\otimes End_{\mathcal{D}}(D)))\simeq Cau(\mathcal{C}\otimes \mathcal{D}).$$ As $\mathds{k}$ is perfect, corollary 18 of \cite{ERZ} implies that $End_{\mathcal{C}}(C)\otimes End_{\mathcal{D}}(D)$ is a finite semisimple algebra. Thus, the left hand-side is the category of finite dimensional modules over a finite semisimple algebra. In particular, we have that $Cau(\mathcal{C}\otimes \mathcal{D})$ is a finite semisimple category. This proves the claim, and concludes the proof.
\end{proof}

We now recall two properties of the Deligne tensor product that will have analogues for the 2-Deligne tensor product. The first result gives us a partial description of the $Hom$-spaces in a Deligne tensor product (for instance, see \cite{EGNO}).

\begin{Lemma}\label{lem:tensorhom}
Given $\mathcal{C}$ and $\mathcal{D}$ two finite categories, for any $C_1, C_2$ in $\mathcal{C}$ and $D_1, D_2$ in $\mathcal{D}$, we have: $$Hom_{\mathcal{C}\boxtimes \mathcal{D}}(C_1\boxtimes D_1, C_2\boxtimes D_2)\cong Hom_{\mathcal{C}}(C_1, C_2)\otimes Hom_{\mathcal{C}}(D_1, D_2).$$
\end{Lemma}

The second tells us how the Deligne tensor product behaves if we input monoidal finite tensor categories, and then finite tensor categories.

\begin{Proposition}\label{prop:1DeligneMon}
Let $\mathcal{C}$ and $\mathcal{D}$ be finite monoidal categories. Then, the finite category $\mathcal{C}\boxtimes\mathcal{D}$, and the 2-universal functor $\boxtimes: \mathcal{C}\times \mathcal{D}\rightarrow \mathcal{C}\boxtimes \mathcal{D}$ admit canonical monoidal structures. Further, $\boxtimes$ is 2-universal as a monoidal bilinear bifunctor right-exact in both variables.
\end{Proposition}

\begin{proof}
It is known that the Deligne tensor product of two finite monoidal categories is monoidal, and that the functor $\boxtimes$ is indeed monoidal (see \cite{Del} or \cite{DSPS13}). Using the 2-universal property of the 2-Deligne tensor product, the remaining part of the statement follows.
\end{proof}

Let us also recall the following result, which is proposition 5.17 of \cite{Del}.

\begin{Proposition}\label{prop:1Delignerigid}
Assume $\mathds{k}$ is perfect. Let $\mathcal{C}$ and $\mathcal{D}$ be finite tensor categories. Then, the canonical monoidal structure on $\mathcal{C}\boxtimes\mathcal{D}$ is rigid, i.e. $\mathcal{C}\boxtimes\mathcal{D}$ is a finite tensor category.
\end{Proposition}

\begin{Remark}
As is argued in subsection 2.2.3 of \cite{DSPS13}, the hypothesis that $\mathds{k}$ is perfect cannot be removed in corollary \ref{prop:1Delignerigid}. This explains why we will focus our attention on perfect fields in what follows.
\end{Remark}

\section{A Review of Compact Semisimple 2-Categories}\label{sec:recollections}

As we have seen above, one important feature of finite semisimple categories is that they are Cauchy complete, i.e. they admit direct sums and idempotents split. It is therefore desirable that any categorification of this notion should be suitably Cauchy complete. Defining directs sums for objects in a 2-category, by which we mean a weak 2-category or bicategory, is not difficult. However, categorifying the notion of idempotents and their splittings requires is not nearly as straightforward. One possible way of achieving this is using separable monads and their splittings as explained in \cite{DR}. However, we prefer the notions of 2-condensation monads and 2-condensations introduced in \cite{GJF} (see subsection 1.1 of \cite{D1} for spelled out versions of these two definitions), as they enjoy better categorical properties. Note that this preference is essentially a matter of taste in the context of semisimple 2-categories as can be seen from theorem 3.3.3 of \cite{GJF}. Now, let $\mathds{k}$ be a field. As it will be extensively used in the next section, we begin by spelling out the 3-universal property of the Cauchy completion of $\mathds{k}$-linear 2-categories in the spirit of definition 1.2.1 of \cite{D1}.

\begin{Definition}
Let $\mathfrak{A}$ be a $\mathds{k}$-linear 2-category which is locally Cauchy complete, i.e. its $Hom$-categories are Cauchy complete. A Cauchy completion is a 3-functor $\kappa:\mathfrak{A}\rightarrow Cau(\mathfrak{A})$ to a locally Cauchy complete and Cauchy complete 2-category $Cau(\mathfrak{A})$ satisfying the following 3-universal property:
\begin{enumerate}

\item For any locally Cauchy complete and Cauchy complete $\mathds{k}$-linear 2-category $\mathfrak{E}$ and every linear 2-functor $F:\mathfrak{A}\rightarrow \mathfrak{E}$, there exists a linear 2-functor $F':Cau(\mathfrak{A})\rightarrow\mathfrak{E}$ and a 2-natural equivalence $u:F'\circ\kappa\Rightarrow F$.

\item For every linear 2-functors $G,H:Cau(\mathfrak{A})\rightarrow \mathfrak{E}$ and 2-natural transformation $t:G\circ\kappa\Rightarrow H\circ\kappa$, there exists a 2-natural transformation $t':G\Rightarrow H$ and an invertible modification $\lambda:t'\circ \kappa\Rrightarrow t$.

\item Furthermore, for any 2-natural transformations $r,s:G\Rightarrow H$ and modification $\lambda:r \circ\kappa\Rightarrow s \circ\kappa$, there exists a unique modification $\lambda':r \Rrightarrow s$ such that $\lambda'\circ \kappa=\lambda$.
\end{enumerate}
\end{Definition}

We move one to recall the definition of a semisimple 2-category over $\mathds{k}$ from \cite{D5} (see \cite{DR} for the original definition).

\begin{Definition}
A semisimple 2-category is a $\mathds{k}$-linear 2-category, which has right and left adjoints for 1-morphisms, is Cauchy complete, and whose $Hom$-categories are semisimple abelian, i.e. every object splits as a finite direct sum of simple objects. 
\end{Definition}

An object $C$ of a semisimple 2-category $\mathfrak{C}$ is called simple if $Id_C$ is a simple object of finite semisimple tensor category $End_{\mathfrak{C}}(C)$. Over an algebraically closed field of characteristic zero, the following definition was introduced in \cite{DR} to categorify the notion of a finite semisimple 1-category. It was generalized in \cite{D5} to arbitrary fields.

\begin{Definition}
A semisimple 2-category is called finite if it is locally finite, i.e. its $Hom$-categories are finite, and it has only finitely many equivalence classes of simple objects.
\end{Definition}

However, as shown in theorem 2.1.2 of \cite{D5}, finite semisimple 2-categories do not exist over fields that are not separably closed. As we wish to work over as general a field as possible, we need to recall a generalization of the notion of finite semisimple 2-category. Namely, following \cite{D5}, we say that two simple objects of $\mathfrak{C}$ are connected if there exists a non-zero 1-morphism between them. As observed in \cite{DR}, being connected defines an equivalence relation on the set of simple objects of $\mathfrak{C}$, and we write $\pi_0(\mathfrak{C})$ for the quotient. We are now ready to recall the definition of a compact semisimple 2-category from \cite{D5}.

\begin{Definition}
A semisimple 2-category $\mathfrak{C}$ is called compact if it is locally finite, and $\pi_0(\mathfrak{C})$ is finite.
\end{Definition}

Given a finite semisimple tensor category $\mathcal{C}$, theorem 1.3.4 of \cite{D5} proves that $Cau(\mathrm{B}\mathcal{C})\simeq \mathbf{Mod}(\mathcal{C})$ is a compact semisimple 2-category, which we can equivalently think of as the 2-category of separable right $\mathcal{C}$-module categories (see definition 2.5.1 of \cite{DSPS13}). Moreover, theorem 1.3.6 of \cite{D5} characterises compact semisimple 2-categories completely, and provides further justification as to why the last definition is the appropriate categorification of the notion of a finite semisimple 1-category over an arbitrary field. We give the statement below for the reader's convenience.

\begin{Theorem}\label{thm:generator}
Let $\mathfrak{C}$ be a compact semisimple 2-category. Then, there exists an object $C$ of $\mathfrak{C}$ such that the inclusion $\mathrm{B}End_{\mathfrak{C}}(C)\hookrightarrow \mathfrak{C}$ is a Cauchy completion.
\end{Theorem}

\noindent In particular, every compact semisimple 2-category $\mathfrak{C}$ is equivalent to $\mathbf{Mod}(\mathcal{C})$ from some finite semisimple tensor category $\mathcal{C}$.

\begin{Remark}
Let us also mention that corollary 2.2.3 of \cite{D5} shows that, over an algebraically closed field or a real closed field, every compact semisimple 2-category is finite.
\end{Remark}

Finally, we will also use the following variant of the notion of compact semisimple 2-category introduced in \cite{D5}.

\begin{Definition}
A compact semisimple 2-category $\mathfrak{C}$ is locally separable if for every object $C$ of $\mathfrak{C}$, $End_{\mathfrak{C}}(C)$ is a separable tensor category (see definition 2.5.8 of \cite{DSPS13}).
\end{Definition}

\section{The 2-Deligne Tensor Product}\label{sec:2Deligne}

Currently, there is not adequate definition of finite 2-category available. In particular, we cannot mimick definition \ref{def:1Deligne} in order to obtain a definition of the 2-Deligne tensor product. However, we can adapt the second construction given in section \ref{sec:1Deligne}. Namely, we begin by defining a completed tensor product of $\mathds{k}$-linear 2-categories over an arbitrary field $\mathds{k}$. Then, we show that the completed tensor product of two compact semisimple 2-categories is compact semisimple 2-category, this proves that the 2-Deligne tensor product of two compact semisimple 2-categories exists.

\begin{Definition}\label{def:2completedtensor}
Given $\mathfrak{A}$ and $\mathfrak{B}$ two $\mathds{k}$-linear 2-categories. If it exists, we write $\widehat{\otimes}: \mathfrak{A}\times\mathfrak{B}\rightarrow \mathfrak{A}\widehat{\otimes}\mathfrak{B}$ for the bilinear 2-functor to a locally Cauchy complete and Cauchy complete $\mathds{k}$-linear 2-category satisfying the following 3-universal property:
\begin{enumerate}
\item For any Cauchy complete and locally Cauchy complete $\mathds{k}$-linear 2-category $\mathfrak{E}$, and every bilinear 2-functor $F:\mathfrak{A}\times\mathfrak{B}\rightarrow \mathfrak{E}$, there exists a linear 2-functor $F':\mathfrak{A}\widehat{\otimes}\mathfrak{B}\rightarrow\mathfrak{E}$ and a 2-natural equivalence $u:F'\circ\widehat{\otimes}\Rightarrow F$.

\item For every 2-functors $G,H:\mathfrak{A}\widehat{\otimes}\mathfrak{B}\rightarrow \mathfrak{E}$, and 2-natural transformation $t:G\circ\widehat{\otimes}\Rightarrow H\circ \widehat{\otimes}$, there exists a 2-natural equivalence $t':G\Rightarrow H$ together with an invertible modification $\gamma:t' \circ \widehat{\otimes} \Rrightarrow t$.

\item Furthermore, for any 2-natural transformation $r,s:G\Rightarrow H$ and modification $\delta:r\circ \widehat{\otimes}\Rrightarrow s\circ\widehat{\otimes}$, there exists a unique invertible modification $\delta:r\Rrightarrow s$ such that $\delta'\circ\widehat{\otimes}=\delta$.
\end{enumerate}
\end{Definition}

\begin{Remark}\label{rem:2Delignereformulation}
Definition \ref{def:2completedtensor} can be reformulated more concisely as saying that precomposition with $\widehat{\otimes}$ induces an equivalence of 2-categories $$Hom(\mathfrak{A}\widehat{\otimes}\mathfrak{B},\mathfrak{E})\simeq Hom_{bil}(\mathfrak{A}\times\mathfrak{B},\mathfrak{E}),$$ from the 2-category of linear 2-functors $\mathfrak{A}\widehat{\otimes}\mathfrak{B}\rightarrow\mathfrak{E}$ to the 2-category of bilinear 2-functors $\mathfrak{A}\times\mathfrak{B}\rightarrow\mathfrak{E}$.
\end{Remark}

We now wish to show that the 3-universal bilinear 2-functor $\widehat{\otimes}: \mathfrak{A}\times\mathfrak{B}\rightarrow \mathfrak{A}\widehat{\otimes}\mathfrak{B}$ of definition \ref{def:2completedtensor} exists. We begin with the following lemma.

\begin{Lemma}\label{lem:locallycompletedtensorproduct}
Given $\mathfrak{A}$, and $\mathfrak{B}$ two $\mathds{k}$-linear 2-categories. There is a 3-universal bilinear 2-functor $\beta:\mathfrak{A}\times \mathfrak{B}\rightarrow \mathfrak{A}\widetilde{\otimes}\mathfrak{B}$, i.e. for every $\mathds{k}$-linear and locally Cauchy complete 2-category $\mathfrak{E}$, precomposition with $\beta$ induces an equivalence of 2-categories
$$Hom(\mathfrak{A}\widetilde{\otimes} \mathfrak{B},\mathfrak{E})\simeq Hom_{bil}(\mathfrak{A}\times \mathfrak{B},\mathfrak{E})$$ from the 2-category of linear 2-functors $\mathfrak{A}\widetilde{\otimes} \mathfrak{B}\rightarrow\mathfrak{E}$ to the 2-category of bilinear 2-functors $\mathfrak{A}\times \mathfrak{B}\rightarrow\mathfrak{E}$.
\end{Lemma}

\begin{proof}
Let $\mathfrak{A}\otimes\mathfrak{B}$ be the $\mathds{k}$-linear 2-category with objects given by the pairs $(A,B)$ with $A$ an object of $\mathfrak{A}$, and $B$ an object of $\mathfrak{B}$, and with $Hom$-categories from $(A_1,B_1)$ to $(A_2,B_2)$ given by $$Hom_{\mathfrak{A}\otimes\mathfrak{B}}((A_1,B_1),(A_2,B_2)):= Hom_{\mathfrak{A}}(A_1,A_2)\otimes Hom_{\mathfrak{B}}(B_1,B_2).$$ The 2-universal property of $\otimes$ explained in lemma \ref{lem:tensorproductcategories} shows that this does define a 2-category, and that precomposition with the canonical bilinear 2-functor $\mathfrak{A}\times \mathfrak{B}\rightarrow \mathfrak{A}\otimes \mathfrak{B}$ induces an equivalence of linear 2-categories $$Hom(\mathfrak{A}\otimes \mathfrak{B},\mathfrak{F})\simeq Hom_{bil}(\mathfrak{A}\times \mathfrak{B},\mathfrak{F})$$ for every linear 2-category $\mathfrak{F}$.

Now, let $\mathfrak{A}\widetilde{\otimes}\mathfrak{B}$ be the local Cauchy completion of $\mathfrak{A}\otimes\mathfrak{B}$, i.e. the 2-category with objects given by the pairs $(A,B)$, and with $Hom$-categories given by the Cauchy completions of $Hom_{\mathfrak{A}\otimes\mathfrak{B}}((A_1,B_1),(A_2,B_2))$. Using the 2-universal property of the Cauchy completion of $\mathds{k}$-linear categories of \cite{D1} recalled in lemma \ref{lem:1CauchyCompletion} above, we see that this does define a 2-category, and that precomposition with the canonical inclusion $\mathfrak{A}\otimes\mathfrak{B}\rightarrow\mathfrak{A}\widetilde{\otimes}\mathfrak{B}$ induces an equivalence of 2-categories 
$$Hom(\mathfrak{A}\widetilde{\otimes} \mathfrak{B},\mathfrak{E})\simeq Hom_{bil}(\mathfrak{A}\otimes \mathfrak{B},\mathfrak{E})$$ for every locally Cauchy complete linear 2-category $\mathfrak{E}$.
\end{proof}

\begin{Proposition}\label{prop:existencecompletedtensor}
For any $\mathds{k}$-linear 2-categories $\mathfrak{A}$ and $\mathfrak{B}$, the 3-universal bilinear 2-functor $\widehat{\otimes}: \mathfrak{A}\times\mathfrak{B}\rightarrow \mathfrak{A}\widehat{\otimes}\mathfrak{B}$ of definition \ref{def:2completedtensor} exists.
\end{Proposition}
\begin{proof}
Using the notations of lemma \ref{lem:locallycompletedtensorproduct}, we set $\mathfrak{A}\widehat{\otimes}\mathfrak{B}$ to be the Cauchy completion of the 2-category $\mathfrak{A}\widetilde{\otimes}\mathfrak{B}$, and we let $\widehat{\otimes}:\mathfrak{A}\times\mathfrak{B}\rightarrow \mathfrak{A}\widehat{\otimes}\mathfrak{B}$ be equal to the composite $$\mathfrak{A}\times\mathfrak{B}\xrightarrow{\beta}\mathfrak{A}\widetilde{\otimes}\mathfrak{B}\xrightarrow{\kappa}Cau(\mathfrak{A}\widetilde{\otimes}\mathfrak{B})=\mathfrak{A}\widehat{\otimes}\mathfrak{B}.$$ Now, for any locally Cauchy complete and Cauchy complete $\mathds{k}$-linear 2-category $\mathfrak{E}$, precomposition by $\widehat{\otimes}$ induces an equivalence of 2-categories $$Hom(\mathfrak{A}\widehat{\otimes}\mathfrak{B},\mathfrak{E})\simeq Hom(\mathfrak{A}\widetilde{\otimes}\mathfrak{B}, \mathfrak{E})\simeq Hom_{bil}(\mathfrak{A}\times\mathfrak{B}, \mathfrak{E}).$$ This proves the result.
\end{proof}

\begin{Remark}
Thanks to its definition using a 3-universal property, the bilinear 2-functor $\widehat{\otimes}: \mathfrak{A}\times\mathfrak{B}\rightarrow \mathfrak{A}\widehat{\otimes}\mathfrak{B}$ is unique in the sense that the 3-category formed by such 2-functors is a contractible 3-groupoid.
\end{Remark}

From now on, we take $\mathds{k}$ to be a perfect field. By analogy with the decategorified situation recalled in section \ref{sec:1Deligne}, we make the following definition.

\begin{Definition}\label{def:2Deligne}
Given two compact semisimple 2-categories $\mathfrak{C}$ and $\mathfrak{D}$ over a perfect field $\mathds{k}$, their 2-Deligne tensor product, if it exists, is the 3-universal bilinear 2-functor $\boxdot:\mathfrak{C}\times\mathfrak{D}\rightarrow \mathfrak{C}\boxdot\mathfrak{D}$ to a compact semisimple 2-category.
\end{Definition}

Observe that, if we can prove that the completed tensor product $\mathfrak{C}\widehat{\otimes}\mathfrak{D}$ is a compact semisimple 2-category, we immediately get that the 2-Deligne tensor product of $\mathfrak{C}$ and $\mathfrak{D}$ exists.

\begin{Theorem}\label{thm:2Deligne}
Over a perfect field, the 2-Deligne tensor product $\boxdot:\mathfrak{C}\times\mathfrak{D}\rightarrow\mathfrak{C}\boxdot\mathfrak{D}$ of two compact semisimple 2-categories $\mathfrak{C}$ and $\mathfrak{D}$ exists.
\end{Theorem}

\begin{proof}
It is enough to prove that $\mathfrak{C}\widehat{\otimes}\mathfrak{D}=Cau(\mathfrak{C}\widetilde{\otimes}\mathfrak{D})$ is a compact semisimple 2-category, as then, $\widehat{\otimes}: \mathfrak{C}\times\mathfrak{D}\rightarrow \mathfrak{C}\widehat{\otimes}\mathfrak{D}$ has the desired 3-universal property. To this end, let $C$ be a generator of $\mathfrak{C}$, and $D$ be a generator of $\mathfrak{D}$, meaning that the canonical inclusions $\kappa_C:\mathrm{B}End_{\mathfrak{C}}(C)\rightarrow\mathfrak{C}$ and $\kappa_D:\mathrm{B}End_{\mathfrak{D}}(D)\rightarrow \mathfrak{D}$ are Cauchy completions. (These exist by theorem \ref{thm:generator}.) Note that the linear 2-functor $$\kappa_C\widetilde{\otimes} \kappa_D:\mathrm{B}End_{\mathfrak{C}}(C)\widetilde{\otimes} \mathrm{B}End_{\mathfrak{D}}(D)\rightarrow \mathfrak{C}\widetilde{\otimes} \mathfrak{D}$$ is fully faithful by construction. We claim that every object of $\mathfrak{C}\widetilde{\otimes} \mathfrak{D}$ is the splitting of a 2-condensation monad supported on the image $C\widetilde{\otimes} D$ of $\kappa_C\widetilde{\otimes} \kappa_D$. To see this, observe that as $C$ is a generator for $\mathfrak{C}$, every simple object of $\mathfrak{C}$ is the splitting of a 2-condensation monad supported on $C$. Now, if $(C,C_1,f_1,g_1,\phi_1,\gamma_1)$ and $(C,C_2,f_2,g_2,\phi_2,\gamma_2)$ are two 2-condensations in $\mathfrak{C}$ (in the notation of subsection 1.1 of \cite{D1}), note that $$\Big(C,C_1\boxplus C_1, \begin{pmatrix}f_1\\f_2\end{pmatrix}, \begin{pmatrix}g_1& g_2\end{pmatrix}, \begin{pmatrix}\phi_1 & 0\\0 & \phi_2\end{pmatrix},\begin{pmatrix}\gamma_1 & 0\\0 & \gamma_2\end{pmatrix}\Big)$$ is a 2-condensation. In particular, $C_1\boxplus C_2$ is the splitting of a 2-condensation monad supported on $C$. As every object of $\mathfrak{C}$ is a finite direct sum of simple objects by proposition 1.1.7 of \cite{D5}, we find that every object $C'$ of $\mathfrak{C}$ is the splitting of a condensation monad supported on $C$. Likewise, every object $D'$ of $\mathfrak{D}$ is the splitting of a 2-condensation monad supported on $D$. Thus, we find that $(C', D')$ is the splitting of a 2-condensation monad supported on $(C, D)$. Now, using that splittings of 2-condensation monads are preserved by all 2-functors, we get that $C'\widetilde{\otimes} D'$ is the splitting of a 2-condensation monad supported on $C\widetilde{\otimes} D$ as desired.

Using lemma A.2.4 of \cite{D1}, we see that $Cau(\kappa_C\widetilde{\otimes} \kappa_D)$ is an equivalence of 2-categories. Moreover, by construction, we have $$BEnd_{\mathfrak{C}}(C)\widetilde{\otimes} BEnd_{\mathfrak{D}}(D)\simeq B(End_{\mathfrak{C}}(C)\boxtimes End_{\mathfrak{D}}(D)).$$ But $End_{\mathfrak{C}}(C)\boxtimes End_{\mathfrak{D}}(D)$ is a finite semisimple tensor category by corollary \ref{prop:1Delignerigid}. This completes the proof.
\end{proof}

\begin{Corollary}\label{cor:2DeligneFSS2C}
Let $\mathcal{C}$ and $\mathcal{D}$ be two finite semisimple tensor categories. Then, we have $$\mathbf{Mod}(\mathcal{C})\boxdot \mathbf{Mod}(\mathcal{D})\simeq \mathbf{Mod}(\mathcal{C}\boxtimes \mathcal{D}).$$
\end{Corollary}

\begin{Corollary}
Let $\mathfrak{C}$ and $\mathfrak{D}$ be two locally separable compact semisimple 2-categories. Then, $\mathfrak{C}\boxdot \mathfrak{D}$ is locally separable.
\end{Corollary}
\begin{proof}
By theorem 1.4.7 of \cite{D5}, there exists separable tensor categories $\mathcal{C}$ and $\mathcal{D}$ such that $\mathfrak{C}\simeq\mathbf{Mod}(\mathcal{C})$ and $\mathfrak{D}\simeq\mathbf{Mod}(\mathcal{D})$. By corollary \ref{cor:2DeligneFSS2C}, we find that $\mathfrak{C}\boxdot\mathfrak{D}\simeq \mathbf{Mod}(\mathcal{C}\boxdot\mathcal{D})$. The result now follows from corollary 2.5.11 of \cite{DSPS13}.
\end{proof}

\begin{Remark}
Observe that under the 3-equivalence of theorem 3.1.4 of \cite{D5}, the 1-Deligne tensor product $\mathcal{C}\boxtimes\mathcal{D}$ of two finite semisimple tensor categories is sent to the 2-Deligne tensor product of their compact semisimple categories of separable modules $\mathbf{Mod}(\mathcal{C})\boxdot \mathbf{Mod}(\mathcal{D})$. Thus, taking for granted that the 3-category $\mathbf{TC}^{rsep}$ of finite semisimple tensor categories and right separable bimodules admits a symmetric monoidal structure given by the Deligne tensor product, then $\mathbf{CSS2C}$ the 3-category of compact semisimple 2-categories admits a symmetric monoidal structure given by the 2-Deligne tensor product. In particular, as every separable tensor category is a fully dualizable object $\mathbf{TC}^{rsep}$ by results of \cite{DSPS13}, every locally separable compact semisimple 2-category is a fully dualizable object $\mathbf{CSS2C}$.
\end{Remark}

\section{Properties of the 2-Deligne Tensor Product}\label{sec:properties}

Now that we have shown that the 2-Deligne tensor product of two compact semisimple 2-categories is a compact semisimple 2-category, we can examine its properties. Throughout, we work over a perfect field $\mathds{k}$ unless otherwise specified.

\begin{Proposition}\label{prop:2DeligneHom}
Let $\mathfrak{C}$ and $\mathfrak{D}$ be two compact semisimple 2-categories. For any objects $C_1, C_2$ in $\mathfrak{C}$ and $D_1, D_2$ in $\mathfrak{D}$, there is an equivalence of categories: $$Hom_{\mathfrak{C}\boxdot \mathfrak{D}}(C_1\boxdot D_1, C_2\boxdot D_2)\simeq Hom_{\mathfrak{C}}(C_1, C_2)\boxtimes Hom_{\mathfrak{C}}(D_1, D_2).$$ Moreover, this equivalence commutes with composition up to natural isomorphism.
\end{Proposition}
\begin{proof}
By theorem 1.3.6 of \cite{D5} recalled in theorem \ref{thm:generator} above, it is enough to prove that given two finite semisimple tensor categories $\mathcal{C}$ and $\mathcal{D}$, and separable algebras $P_1, P_2\in\mathcal{C}$ and $Q_1, Q_2\in\mathcal{D}$, then $$Bimod_{\mathcal{C}\boxtimes\mathcal{D}}(P_1\boxtimes Q_1, P_2\boxtimes Q_2)\simeq Bimod_{\mathcal{C}}(P_1, P_2)\boxtimes Bimod_{\mathcal{D}}(Q_1, Q_2).$$ Observe that there is a natural functor $F$, as in the diagram below, given by taking the Deligne tensor product of two bimodules. As $F$ is bilinear, it can be factored through the Deligne tensor product up to natural isomorphism, i.e. there exist a functor $F'$ making the diagram $$\begin{tikzcd}
{Bimod_{\mathcal{C}}(P_1, P_2)\times Bimod_{\mathcal{D}}(Q_1, Q_2)} \arrow[r, "F"] \arrow[d, "\boxtimes"'] &                                                                                          {Bimod_{\mathcal{C}\boxtimes\mathcal{D}}(P_1\boxtimes Q_1, P_2\boxtimes Q_2)} \\
                                                                                                             {Bimod_{\mathcal{C}}(P_1, P_2)\boxtimes Bimod_{\mathcal{D}}(Q_1, Q_2).} \arrow[ru, dotted, "F'"'] &                                                                              
\end{tikzcd}$$ commute up to natural isomorphism.

Let us prove that the functor $F'$ is fully faithful. Let $M$ be an object of $ Bimod_{\mathcal{C}}(P_1, P_2)$, $N$ an object of $Bimod_{\mathcal{D}}(Q_1, Q_2)$, and $C\boxtimes D$ an object of $\mathcal{C}\boxtimes\mathcal{D}$. There are natural isomorphisms:
\begin{align*}
    Hom_{P_1\boxtimes Q_1 , P_2\boxtimes Q_2}&((P_1\boxtimes Q_1)\otimes (C\boxtimes D) \otimes (P_2\boxtimes Q_2), M\boxtimes N)\\
    &\cong Hom_{\mathcal{C}\boxtimes \mathcal{D}}(C\boxtimes D, M\boxtimes N)\\
    &\cong Hom_{\mathcal{C}}(C, M)\otimes Hom_{\mathcal{D}}(D, N)\\
    &\cong Hom_{P_1 , P_2}(P_1\otimes C\otimes P_2, M)\otimes Hom_{Q_1 , Q_2}(Q_1\otimes D\otimes Q_2, N).
\end{align*}
As every object of $Bimod_{\mathcal{C}}(P_1, P_2)\boxtimes Bimod_{\mathcal{D}}(Q_1, Q_2)$ is the splitting of an idempotent on an object of the form $(P_1\otimes C\otimes P_2)\boxtimes (Q_1\otimes D\otimes Q_2)$ for some objects $C$ in $\mathcal{C}$ and $D$ in $\mathcal{D}$, we get that the functor $F'$ is fully faithful.

Now, observe that the objects $(P_1\boxtimes Q_1)\otimes (C\boxtimes D) \otimes (P_2\boxtimes Q_2)$ with $C\boxtimes D$ varying form a set of projective generators for the finite semisimple category $Bimod_{\mathcal{C}\boxtimes\mathcal{D}}(P_1\boxtimes Q_1, P_2\boxtimes Q_2)$. Given that we have shown above that $F'$ is fully faithful, this implies that $F'$ is essentially surjective.

Finally, let us comment on the naturality of this equivalence. Given $P_1, P_2, P_3$ separable algebras in $\mathcal{C}$, and $Q_1,Q_2,Q_3$ separable algebras in $\mathcal{D}$, for every bimodules $A_1$ in $Bimod(P_1, P_2)$, $A_2$ in  $Bimod(P_2, P_3)$, $B_1$ in $Bimod(Q_1,Q_2)$, and $B_2$ in $Bimod(Q_2,Q_3)$, using the 2-universal property of the Deligne tensor product, it is enough to show that there is a natural isomorphism $$(A_1\otimes_{P_2}A_2)\boxtimes (B_1\otimes_{Q_2}B_2)\cong (A_1\boxtimes B_1)\otimes_{P_2\boxtimes Q_2}(A_2\boxtimes B_2).$$
This follows by comparing the universal properties.
\end{proof}

\begin{Lemma}
Let us fix a compact semisimple 2-category $\mathfrak{D}$. Then, for every compact semisimple 2-categories $\mathfrak{C}$, and $\mathfrak{E}$, there is a natural equivalence of 2-categories
$$Hom(\mathfrak{C}\boxdot \mathfrak{D},\mathfrak{E})\simeq Hom(\mathfrak{C}, Hom(\mathfrak{D},\mathfrak{E})).$$
\end{Lemma}
\begin{proof}
We have the following natural equivalences of 2-categories:
$$Hom(\mathfrak{C}\boxdot \mathfrak{D},\mathfrak{E})\simeq Hom_{bil}(\mathfrak{C}\times \mathfrak{D},\mathfrak{E}) \simeq Hom(\mathfrak{C}, Hom(\mathfrak{D},\mathfrak{E})),$$ where the first equivalence is given by theorem \ref{rem:2Delignereformulation}, and the second holds by inspection.
\end{proof}

One may also ask how the connected components of a compact semisimple 2-category interact with the 2-Deligne tensor product. Over an algebraically closed field, this is the same thing as asking what are the connected components of the 2-Deligne tensor product of two finite semisimple 2-categories. In this case, we simply get the product of the connected components of the finite semisimple 2-categories.

\begin{Lemma}\label{lem:connected2Deligne}
Let $\mathds{k}$ be an algebraically closed field. Let $\mathfrak{C}$, $\mathfrak{D}$ be two finite semisimple 2-categories. Then, we have $$\pi_0(\mathfrak{C}\boxdot\mathfrak{D}) \cong \pi_0(\mathfrak{C})\times \pi_0(\mathfrak{D}).$$
\end{Lemma}
\begin{proof}
Observe that any finite semisimple 2-category can be decomposed into a direct sum of connected finite semisimple 2-categories (i.e. $\pi_0$ is a singleton). Therefore, it is enough to prove that if $\mathfrak{C}$ and $\mathfrak{D}$ are connected, then so $\mathfrak{C}\boxdot\mathfrak{D}$. But, a finite semisimple 2-category is connected if and only if it is equivalent to the 2-category of separable module categories over a fusion category. Thus, there exists fusion categories $\mathcal{C}$ and $\mathcal{D}$ such that $\mathfrak{C}\simeq\mathbf{Mod}(\mathcal{C})$ and $\mathfrak{D}\simeq\mathbf{Mod}(\mathcal{D})$. By corollary \ref{cor:2DeligneFSS2C}, we have $$\mathfrak{C}\boxdot\mathfrak{D}\simeq \mathbf{Mod}(\mathcal{C}\boxtimes\mathcal{D}).$$ The result now follows from lemma \ref{lem:tensorhom}, which implies that $\mathcal{C}\boxtimes\mathcal{D}$ is a fusion category.
\end{proof}

\begin{Example}\label{ex:VectR}
Let us work over $\mathbb{R}$. Then, we have \begin{align*}\mathbf{Mod}(\mathbf{Vect}_{\mathbb{C}})&\boxdot \mathbf{Mod}(\mathbf{Vect}_{\mathbb{C}})\simeq \mathbf{Mod}(\mathbf{Vect}_{\mathbb{C}}\boxtimes \mathbf{Vect}_{\mathbb{C}})\\ &\simeq \mathbf{Mod}(\mathbf{Vect}_{\mathbb{C}}\boxplus \mathbf{Vect}_{\mathbb{C}})\simeq \mathbf{Mod}(\mathbf{Vect}_{\mathbb{C}})\boxplus \mathbf{Mod}(\mathbf{Vect}_{\mathbb{C}}).\end{align*} In particular, lemma \ref{lem:connected2Deligne} fails over non-algebraically closed fields. Moreover, note that $\mathbf{Vect}_{\mathbb{C}}\boxdot \mathbf{Vect}_{\mathbb{C}}$ is the direct sum of the two distinct simple objects of $\mathbf{Mod}(\mathbf{Vect}_{\mathbb{C}})\boxplus \mathbf{Mod}(\mathbf{Vect}_{\mathbb{C}})$. Thus, the 2-functor $\boxdot$ may send the product of two simple objects to an object that is not simple. (Note however that this can not happen over an algebraically closed field as, in this case, the Deligne tensor product of two fusion categories is a fusion category by lemma \ref{lem:tensorhom}.)
\end{Example}

It is useful to know when the functor $\boxdot$ hits every equivalence class of simple objects of the target. Over an algebraically closed field, we give a criterion ensuring this property. Firstly, recall from subsection 4.5 \cite{EGNO} that the Frobenus-Perron dimension of a finite tensor category is an algebraic integer that only depends on the Grothendieck ring. Secondly, observe that given $\mathfrak{C}$ a connected finite semisimple 2-category, we have that $Hom_{\mathfrak{C}}(C,D)$ is non-zero for every pair of simple objects $C,D$. Thus, the endomorphism fusion categories of any two simple objects are Morita equivalent, whence they have the same Frobenius-Perron dimension by corollary 7.16.7 of \cite{EGNO}. This means that the following definition is sensible.

\begin{Definition}
The Frobenius-Perron dimension of a connected finite semisimple 2-category over an algebraically closed field is the Frobenius-Perron dimension of the endomorphism fusion category of any of its simple objects.
\end{Definition}

\begin{Theorem}\label{thm:2Deligneessentiallysurjective}
Let $\mathfrak{C}$ and $\mathfrak{D}$ be connected finite semisimple 2-categories over an algebraically closed field. If $\mathfrak{C}$ and $\mathfrak{D}$ have coprime Frobenius-Perron dimension, then, for any simple object $A$ of $\mathfrak{C}\boxdot\mathfrak{D}$, there exists simple objects objects $C$ in $\mathfrak{C}$ and $D$ in $\mathfrak{D}$ unique up to equivalence such that $C\boxdot D \simeq A$.
\end{Theorem}
\begin{proof}
The theorem follows readily from the following claim: Given $\mathcal{C}$ and $\mathcal{D}$ two fusion categories of coprime Frobenius-Perron dimension. Let $\mathcal{P}$ be an indecomposable right $\mathcal{C}\boxtimes\mathcal{D}$-module category. Then, there exists indecomposable right module categories $\mathcal{M}$ over $\mathcal{C}$ and $\mathcal{N}$ over $\mathcal{D}$ unique up to equivalence such that $\mathcal{P}\simeq \mathcal{M}\boxtimes \mathcal{N}$. Furthermore, if $\mathcal{P}$ is separable, so are $\mathcal{M}$ and $\mathcal{N}$. We prove the first part of this statement using a variant of the proof of proposition 8.55 of \cite{ENO}, which deals with the characteristic zero case.

Given $\mathcal{P}$ an indecomposable right $\mathcal{C}\boxtimes\mathcal{D}$-module category, we can restrict $\mathcal{P}$ to $\mathcal{C}$, and write $\mathcal{P}|_{\mathcal{C}}=\oplus_i\mathcal{M}_i$ for some indecomposable right $\mathcal{C}$-module categories $\mathcal{M}_i$. Then, the action of $\mathcal{D}$ can be viewed as a tensor functor $F:\mathcal{D}\rightarrow \oplus_{i,j}Fun_{\mathcal{C}}(\mathcal{M}_i,\mathcal{M}_j)$. We write $Im(F)$ for the image of $F$, i.e. the multifusion category generated by the summands of the objects in $F(\mathcal{D})$. Now, if $Im(F)$ was decomposable, a decomposition of $Im(F)$ would induce a decomposition of $\mathcal{P}$ as a right $\mathcal{C}\boxtimes\mathcal{D}$-module category. By hypothesis, this is not the case, thus $Im(F)$ is indecomposable. In particular, if we write $Im(F)=\oplus_{i,j}Im(F)_{ij}$ with $Im(F)_{ij}$ a sub-category of $Fun_{\mathcal{C}}(\mathcal{M}_i,\mathcal{M}_j)$, then we find that $Im(F)_{ij}$ is non-zero for every $i$ and $j$.

Now, as $F:\mathcal{D}\rightarrow Im(F)$ is surjective (see definition 1.8.3 of \cite{EGNO}), it follows from theorem 6.2.1 of \cite{EGNO} that the Frobenius-Perron dimension of $Im(F)$ divides that of $\mathcal{D}$. But, for every $i$, the fusion subcategory $Im(F)_{ii}$ of $Im(F)$ is also a fusion subcategory of $Fun_{\mathcal{C}}(\mathcal{M}_i,\mathcal{M}_i)$. Thus, the Frobenius-Perron dimension of $Im(F)_{ii}$ divides the Frobenius-Perron dimensions of $\mathcal{C}$ and $\mathcal{D}$. (This follows from corollary 7.16.7 and theorem 7.17.4 of \cite{EGNO}.) As these algebraic integers are coprime, the Frobenius-Perron dimension of $Im(F)_{ii}$ must be equal to $1$. Whence, $Im(F)_{ii}\simeq \mathbf{Vect}$, and so $Im(F)_{ij}\simeq \mathbf{Vect}$ for every $i$ and $j$, as we are working over an algebraically closed field. This implies that the $\mathcal{M}_i$ must all be equivalent to a single right $\mathcal{C}$-module category $\mathcal{M}$. The functor $F:\mathcal{D}\rightarrow \oplus_{i,j}\mathbf{Vect}=End(\boxplus_{i}\mathbf{Vect})$ then defines a finite semisimple right $\mathcal{D}$-module category $\mathcal{N}$, and we have $\mathcal{P}\simeq \mathcal{M}\boxtimes \mathcal{N}$. Now, it is clear from the construction that $\mathcal{M}$ and $\mathcal{N}$ are unique.

Let us now assume that $\mathcal{P}$ is separable. We wish to show that $\mathcal{M}$ and $\mathcal{N}$ are separable. Let us begin by proving that $\mathcal{M}$ is separable. In order to see this, observe that there is an equivalence of finite tensor categories $$\mathcal{Z}_{\mathcal{D}}(End_{\mathcal{C}}(\mathcal{P}))\simeq End_{\mathcal{C}\boxtimes\mathcal{D}}(\mathcal{P})$$ between the relative Drinfel'd center of $End_{\mathcal{C}}(\mathcal{P})$ over $F:\mathcal{D}\rightarrow End_{\mathcal{C}}(\mathcal{P})$, and the finite tensor category $End_{\mathcal{C}\boxtimes\mathcal{D}}(\mathcal{P})$. But, $End_{\mathcal{C}\boxtimes\mathcal{D}}(\mathcal{P})$ is a fusion category as $\mathcal{P}$ is an indecomposable separable right $\mathcal{C}\boxtimes\mathcal{D}$-module category. Now, the forgetful functor $$\mathcal{Z}_{\mathcal{D}}(End_{\mathcal{C}}(\mathcal{P}))\rightarrow End_{\mathcal{C}}(\mathcal{P})$$ is surjective. Thence, we find that $End_{\mathcal{C}}(\mathcal{P})$ is finite semisimple by theorem 6.1.16 of \cite{EGNO}. Furthermore, for any $i$, there is a surjective tensor functor $End_{\mathcal{C}}(\mathcal{P})\rightarrow End_{\mathcal{C}}(\mathcal{M}_i)$, so we find that $\mathcal{M}_i\simeq \mathcal{M}$ is separable. Finally, in order to prove that $\mathcal{N}$ is separable, one can proceed similarly. Namely, by uniqueness of $\mathcal{M}$ and $\mathcal{N}$, $\mathcal{P}|_{\mathcal{D}}$ decomposes into a direct sum of right $\mathcal{D}$-module categories which are all equivalent to $\mathcal{N}$. The above argument then shows that $\mathcal{N}$ is indeed separable.
\end{proof}

\begin{Remark}
Example \ref{ex:VectR} shows that theorem \ref{thm:2Deligneessentiallysurjective} fails over fields that are not algebraically closed.
\end{Remark}

\begin{Example}\label{ex:Mod(6)}
Let $\mathds{k}$ be an algebraically closed field of characteristic different from $2$ and $3$. Let us consider the finite semisimple 2-category of modules over the separable fusion category $\mathcal{C}=\mathbf{Vect}_{\mathds{Z}/6\mathds{Z}}$. Thanks to our assumption on the characteristic of $\mathds{k}$, all these module categories are separable. Now, by example 7.4.10 of \cite{EGNO} (see also \cite{Nat}), we know that the equivalence classes of simple objects of $\mathbf{Mod}(\mathcal{C})$ are given by pair consisting of a subgroup of $\mathds{Z}/6\mathds{Z}$ and a cohomology class on this subgroup. Denoting by $triv$ the trivial cohomology class, we get the following list:

\begin{center}
\begin{tabular}{ c c c c }
  $(\langle 0 \rangle,triv)$ & $(\langle 3 \rangle,triv)$ & $(\langle 2 \rangle,triv)$ & $(\langle 1 \rangle,triv).$
\end{tabular}
\end{center}

\noindent Now, observe that we have an equivalence $$\mathcal{C}\simeq \mathbf{Vect}_{\mathds{Z}/2\mathds{Z}}\boxtimes \mathbf{Vect}_{\mathbb{Z}/3\mathbb{Z}}$$ of fusion categories. Thus, combining proposition \ref{prop:2DeligneHom}, theorem \ref{thm:2Deligneessentiallysurjective}, we find that the equivalence classes of simple objects of $\mathbf{Mod}(\mathcal{C})$ are described by pair consisting of an equivalence class of simple objects in $\mathbf{Mod}(\mathbf{Vect}_{\mathbb{Z}/2\mathbb{Z}})$ and in $\mathbf{Mod}(\mathbf{Vect}_{\mathbb{Z}/3\mathbb{Z}})$. We derive the following description in the notation of example 2.1.12 of \cite{D1}:
$$(\mathbf{Vect}_{\mathbb{Z}/2\mathbb{Z}},\mathbf{Vect}_{\mathbb{Z}/3\mathbb{Z}})\ \ \ (\mathbf{Vect}_{\mathbb{Z}/3\mathbb{Z}},\mathbf{Vect})\ \ \  (\mathbf{Vect},\mathbf{Vect}_{\mathbb{Z}/2\mathbb{Z}})\ \ \  (\mathbf{Vect},\mathbf{Vect}).$$

\noindent The corresponding irreducible right $\mathcal{C}$-module categories are given by taking the 1-Deligne tensor product of the above pairs. In this case, these can be made very explicit, and we obtain the following list of right $\mathcal{C}$-module categories:

\begin{center}
\begin{tabular}{ c c c c }
 $\mathbf{Vect}_{\mathds{Z}/6\mathds{Z}}$ & $\mathbf{Vect}_{\mathds{Z}/3\mathds{Z}}$ & $\mathbf{Vect}_{\mathds{Z}/2\mathds{Z}}$ & $\mathbf{Vect}.$
\end{tabular}
\end{center}
The module structures arise from the obvious tensor functors.
\end{Example}

\begin{Example}\label{ex:Mod(2+2)}
Let $\mathds{k}$ be an algebraically closed field of characteristic not equal to $2$. We now consider the finite semisimple 2-category of right module categories over $\mathcal{D}=\mathbf{Vect}_{\mathbb{Z}/2\mathbb{Z}\oplus \mathbb{Z}/2\mathbb{Z}}$. Let us write $a=(1,0)$, $b=(0,1)$, and $c=(1,1)$ the non-zero elements of $\mathbb{Z}/2\mathbb{Z}\oplus \mathbb{Z}/2\mathbb{Z}$. The equivalence classes of simple objects of $\mathbf{Mod}(\mathcal{D})$ are given as follows:

\begin{center}
\begin{tabular}{ c c c }
$(\langle 0 \rangle,triv)$ & $(\langle a \rangle,triv)$ & $(\langle b \rangle,triv)$\\
 $(\langle c \rangle,triv)$ & $(\langle a,b \rangle,triv)$ & $(\langle a,b \rangle,\nu).$
 
\end{tabular}
\end{center}

\noindent We have used $\nu$ to denote the non-trivial class in $H^{2}(\mathbb{Z}/2\mathbb{Z}\oplus \mathbb{Z}/2\mathbb{Z},\mathds{k}^*)$. Observe that we have an equivalence $$\mathcal{D}\simeq \mathbf{Vect}_{\mathbb{Z}/2\mathbb{Z}}\boxtimes \mathbf{Vect}_{\mathbb{Z}/2\mathbb{Z}}$$ of fusion categories. By construction, the essential image of the canonical 2-functor $$\boxdot: \mathbf{Mod}(\mathbf{Vect}_{\mathds{Z}/2\mathbb{Z}})\times \mathbf{Mod}(\mathbf{Vect}_{\mathds{Z}/2\mathbb{Z}})\rightarrow \mathbf{Mod}(\mathcal{D})$$ contains exactly the following equivalence classes of simple objects:

\begin{center}
\begin{tabular}{ c c c c }
$(\langle 0 \rangle,triv)$ & $(\langle a \rangle,triv)$ & $(\langle b \rangle,triv)$ & $(\langle a,b \rangle,triv).$
\end{tabular}
\end{center}

\noindent In particular, the image of $\boxdot$ does not contain every equivalence class of simple objects.
\end{Example}

\section{The 2-Deligne Tensor Product of Compact Semi\-simple Tensor 2-Categories}\label{sec:2DeligneMonoidal}

Let $\mathds{k}$ be an arbitrary field. We prove that the completed tensor product $\widehat{\otimes}$ of two monoidal linear 2-categories inherits a monoidal structure. But, before doing so, let us point out the following lemma.

\begin{Lemma}\label{lem:multi2Deligne}
Let $\mathfrak{A}_1$, ..., $\mathfrak{A}_n$ be linear 2-categories. There is a 3-universal multilinear 2-functor, which we denote by $$\widehat{\otimes}:\mathfrak{A}_1\times...\times\mathfrak{A}_n\rightarrow \mathfrak{A}_1\widehat{\otimes}...\widehat{\otimes}\mathfrak{A}_n.$$ In particular, all the different ways of parenthesizing $\mathfrak{A}_1\widehat{\otimes}...\widehat{\otimes}\mathfrak{A}_n$ are canonically equivalent.
\end{Lemma}
\begin{proof}
One can define the functor $\widehat{\otimes}$ inductively by $$\xymatrix{\mathfrak{A}_1\times...\times\mathfrak{A}_n\ar[rr]^-{\widehat{\otimes}\times id\times ...\times id}&&\mathfrak{A}_1\widehat{\otimes}\mathfrak{A}_2\times...\times\mathfrak{A}_n\ar[r]&...\ar[r]^-{\widehat{\otimes}}&\mathfrak{A}_1\widehat{\otimes}...\widehat{\otimes}\mathfrak{A}_n.}$$ The 3-universal property is derived by induction.
\end{proof}

\begin{Theorem} \label{thm:CompletedTensorMon}
Let $\mathfrak{A}$, $\mathfrak{B}$ be monoidal linear 2-categories. Then, the compact semisimple 2-category $\mathfrak{A}\widehat{\otimes}\mathfrak{B}$ admits a canonical monoidal structure such that the bilinear 2-functor $\widehat{\otimes}:\mathfrak{A}\times\mathfrak{B}\rightarrow\mathfrak{A}\widehat{\otimes} \mathfrak{B}$ is monoidal. Further, the bilinear monoidal 2-functor $\widehat{\otimes}$ is the 3-universal bilinear monoidal 2-functor from $\mathfrak{A}\times\mathfrak{B}$ to a locally Cauchy complete and Cauchy complete monoidal linear 2-category.
\end{Theorem}

\begin{Remark}\label{rem:remadjointequivalences}
In the proof below, we will construct many adjoint 2-natural equivalences. In order not complicate notations unnecessarily, we shall only construct one of the two adjoints; The others can be constructed similarly using the 3-universal property of the 2-Deligne tensor product supplied by theorem \ref{thm:2Deligne}.
\end{Remark}

\begin{proof}
Concisely, the idea behind the proof is to observe that the Cartesian product $\mathfrak{A}\times \mathfrak{B}$ is monoidal, and use the 3-universal property of $\widehat{\otimes}$ to transfer this monoidal structure onto $\mathfrak{A}\widehat{\otimes}\mathfrak{B}$. In more details, in the notation of \cite{SP}, write $(\mathfrak{A}, \Box_{\mathfrak{A}}, I_{\mathfrak{A}}, \alpha_{\mathfrak{A}}, l_{\mathfrak{A}}, r_{\mathfrak{A}})$, respectively $(\mathfrak{B}, \Box_{\mathfrak{B}}, I_{\mathfrak{B}}, \alpha_{\mathfrak{B}}, l_{\mathfrak{B}}, r_{\mathfrak{B}})$ for the monoidal structure on $\mathfrak{A}$, respectively $\mathfrak{B}$. Let us consider the following solid arrow diagram:
$$\xymatrix{
\mathfrak{A}\times\mathfrak{B}\times \mathfrak{A}\times\mathfrak{B}\ar[d]_-{\Box_{\mathfrak{A}\times\mathfrak{B}}}\ar[r]^-{\widehat{\otimes}\times\widehat{\otimes}}&(\mathfrak{A}\widehat{\otimes}\mathfrak{B})\times (\mathfrak{A}\widehat{\otimes}\mathfrak{B})\ar[r]^-{\widehat{\otimes}}&\mathfrak{A}\widehat{\otimes}\mathfrak{B}\widehat{\otimes} \mathfrak{A}\widehat{\otimes}\mathfrak{B}\ar@{.>}[d]^-{\Box}\\\mathfrak{A}\times\mathfrak{B}\ar[rr]_-{\widehat{\otimes}}&&\mathfrak{A}\widehat{\otimes}\mathfrak{B},
}$$
where $\Box_{\mathfrak{A}\times\mathfrak{B}}$ is the composite of the canonical permutation 2-functor $$\mathfrak{A}\times \mathfrak{B}\times\mathfrak{A}\times\mathfrak{B}\rightarrow \mathfrak{A}\times \mathfrak{A}\times\mathfrak{B}\times\mathfrak{B},$$ and the 2-functor $\Box_{\mathfrak{A}}\times\Box_{\mathfrak{B}}$. Using proposition \ref{lem:multi2Deligne}, the dotted 2-functor $\Box$ exists, and there is an 2-natural equivalence $\chi$ witnessing the commutativity of the diagram.

\begin{figure}[!htb]
$$\begin{tikzcd}[column sep=1ex]
                                                                                                                            & (\mathfrak{A}\times\mathfrak{B})^{\times 3} \arrow[ld, "\Box_{\mathfrak{A}\times\mathfrak{B}}\times Id"'] \arrow[rdd, near start, "Id\times\Box_{\mathfrak{A}\times\mathfrak{B}}", crossing over] \arrow[rrr, "\widehat{\otimes}^3"] &                                                                                                                            &                                                                                                    & (\mathfrak{A}\widehat{\otimes}\mathfrak{B})^{\times 3} \arrow[ld, "\Box_{\mathfrak{A}\widehat{\otimes}\mathfrak{B}}\times Id"', near end] \arrow[rdd, "Id\times\Box_{\mathfrak{A}\widehat{\otimes}\mathfrak{B}}"] &                                                                                                   \\
(\mathfrak{A}\times\mathfrak{B})^{\times 2} \arrow[rdd, "\Box_{\mathfrak{A}\times\mathfrak{B}}"'] \arrow[rrr,near end, "\widehat{\otimes}^2"',crossing over] &                                                                                                                                                                                                 &                                                                                                                            & (\mathfrak{A}\widehat{\otimes}\mathfrak{B})^{\times 2}  &                                                                                                                                                                          &                                                                                                   \\
                                                                                                                            &                                                                                                                                                                                                 & (\mathfrak{A}\times\mathfrak{B})^{\times 2} \arrow[ld, "\Box_{\mathfrak{A}\times\mathfrak{B}}"'] \arrow[rrr,near start, "\widehat{\otimes}^2"'] &                                                                                                    &                                                                                                                                                                          & (\mathfrak{A}\widehat{\otimes}\mathfrak{B})^{\times 2} \arrow[ld, "\Box_{\mathfrak{A}\widehat{\otimes}\mathfrak{B}}"] \\
                                                                                                                            & \mathfrak{A}\times\mathfrak{B} \arrow[rrr,"\widehat{\otimes}"']                                                                                                                                           &                                                                                                                            &                                                                                                    & \mathfrak{A}\widehat{\otimes}\mathfrak{B}\arrow[from= luu,near start, "\Box_{\mathfrak{A}\widehat{\otimes}\mathfrak{B}}",crossing over].  &                            
\end{tikzcd}$$
\caption{Construction of the associator (Part 1)}
\label{fig:figure1}
\end{figure}
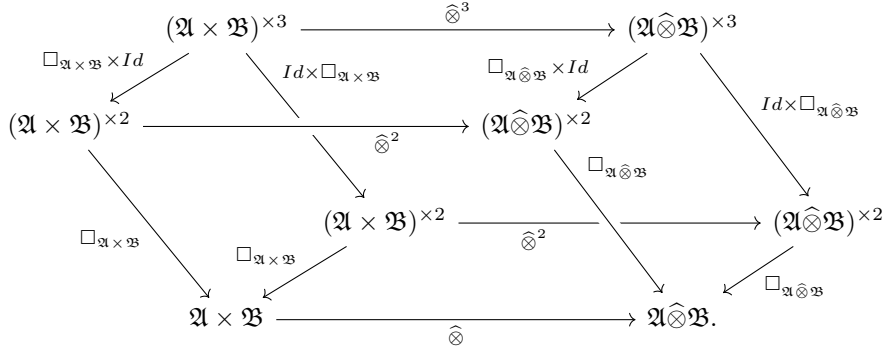

We claim that the 2-functor $\Box_{\mathfrak{A}\widehat{\otimes}\mathfrak{B}}:=\Box\circ\widehat{\otimes}$ provides $\mathfrak{A}\widehat{\otimes}\mathfrak{B}$ with a monoidal structure. Let us show how the associator is constructed. We begin by considering the diagram depicted in figure \ref{fig:figure1}. The left face is filled using the associator of $\mathfrak{A}\times\mathfrak{B}$, and the front, back, top, and bottom faces are filled using $\chi$. Now, we want to use the 3-universal property of the completed tensor product derived in proposition \ref{prop:existencecompletedtensor} to fill in the right face. However, in order for this operation to be licit, we need to partially fill the right face using the diagram depicted in figure \ref{fig:figure2}.

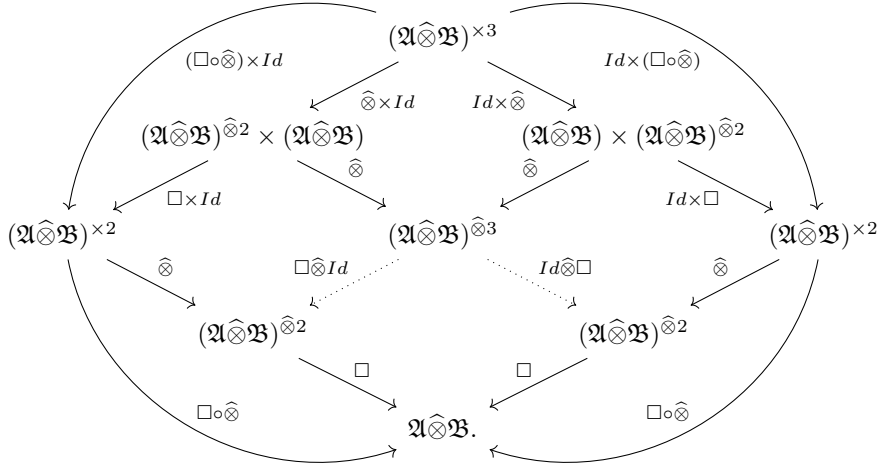
\begin{figure}[!hbt]
$$\adjustbox{max width=\textwidth}{
\begin{tikzcd}[column sep=0]
                                                                                                                   &                                                                                                                                                                                & (\mathfrak{A}\widehat{\otimes}\mathfrak{B})^{\times 3} \arrow[ld, "\widehat{\otimes}\times Id"] \arrow[rd, "Id\times\widehat{\otimes}"'] \arrow[lldd, "(\Box\circ\widehat{\otimes})\times Id", bend right=49] \arrow[rrdd, "Id\times (\Box\circ\widehat{\otimes})"', bend left=49] &                                                                                                                                                                                 &                                                                                                                     \\
                                                                                                                   & (\mathfrak{A}\widehat{\otimes}\mathfrak{B})^{\widehat{\otimes} 2}\times (\mathfrak{A}\widehat{\otimes}\mathfrak{B}) \arrow[rd, "\widehat{\otimes}"] \arrow[ld, "\Box\times Id"] &                                                                                                                                                                                                                                  & (\mathfrak{A}\widehat{\otimes}\mathfrak{B})\times (\mathfrak{A}\widehat{\otimes}\mathfrak{B})^{\widehat{\otimes} 2} \arrow[ld, "\widehat{\otimes}"'] \arrow[rd, "Id\times\Box"'] &                                                                                                                     \\
(\mathfrak{A}\widehat{\otimes}\mathfrak{B})^{\times 2} \arrow[rd, "\widehat{\otimes}"] \arrow[rrdd, "\Box\circ\widehat{\otimes}", bend right=49] &                                                                                                                                                                                & (\mathfrak{A}\widehat{\otimes}\mathfrak{B})^{\widehat{\otimes} 3} \arrow[ld, dotted, "\Box\widehat{\otimes} Id"'] \arrow[rd, dotted, "Id\widehat{\otimes}\Box"]                                                                                                                          &                                                                                                                                                                                 & (\mathfrak{A}\widehat{\otimes}\mathfrak{B})^{\times 2} \arrow[ld, "\widehat{\otimes}"'] \arrow[lldd, "\Box\circ\widehat{\otimes}"', bend left=49] \\
                                                                                                                   & (\mathfrak{A}\widehat{\otimes}\mathfrak{B})^{\widehat{\otimes} 2} \arrow[rd, "\Box"]                                                                                                               &                                                                                                                                                                                                                                  & (\mathfrak{A}\widehat{\otimes}\mathfrak{B})^{\widehat{\otimes} 2} \arrow[ld, "\Box"']                                                                                                               &                                                                                                                     \\
                                                                                                                   &                                                                                                                                                                                & \mathfrak{A}\widehat{\otimes}\mathfrak{B}.                                                                                                                                                                                                  &                                                                                                                                                                                 &                                                                                                                    
\end{tikzcd}
}$$
\caption{Construction of the associator (Part 2)}
\label{fig:figure2}
\end{figure}

There is a canonical 2-natural equivalence filling in the top square of figure \ref{fig:figure2} by lemma \ref{lem:multi2Deligne}. The existence of the dotted arrows, as well as that of 2-natural transformations witnessing the commutativity of the lateral diamonds follow using the 3-universal property of the completed tensor product (see definition \ref{def:2completedtensor}). Observe that the bottom square of the diagram depicted in figure \ref{fig:figure2} is empty.

Gluing together the right face of the diagram represented in figure \ref{fig:figure1} and the exterior face of that given in figure \ref{fig:figure2}, we obtain a diagram on which we can use the 3-universal property of lemma \ref{lem:multi2Deligne}. Doing so, we find a 2-natural equivalence $\alpha_{\mathfrak{A}\widehat{\otimes}\mathfrak{B}}$ that fills the empty square in figure \ref{fig:figure2}, and the the commutativity of the whole polyhedra is witnessed by an invertible modification $\omega$. The modification $\omega$ shows that $\widehat{\otimes}$ and $\chi$ are compatible with the associators. Similarly, one constructs the pentagonator $\pi_{\mathfrak{A}\widehat{\otimes}\mathfrak{B}}$ of $\mathfrak{A}\widehat{\otimes}\mathfrak{B}$ out of the pentagonator of $\mathfrak{A}\times \mathfrak{B}$, simultaneously proving that $\widehat{\otimes}$, $\chi$, and $\omega$ are compatible with the pentagonators. One can then proceed to show that it satisfies the required axioms. Let us briefly comment on the monoidal unit $I_{\mathfrak{A}\widehat{\otimes}\mathfrak{B}}$ of $\mathfrak{A}\widehat{\otimes}\mathfrak{B}$ and on the unitors. In order to defines $I_{\mathfrak{A}\widehat{\otimes}\mathfrak{B}}$, it is convenient to think of objects of a compact semisimple 2-category as linear 2-functors from $\mathrm{B}^2\mathds{k}$, the linear 2-category with one object and one 1-morphism. Consider the following solid arrow diagram: $$\begin{tikzcd}
\mathrm{B}^2\mathds{k}\times \mathrm{B}^2\mathds{k} \arrow[d, "I_{\mathfrak{A}}\times I_{\mathfrak{B}}"'] \arrow[r, "\otimes"] & \mathrm{B}^2\mathds{k} \arrow[d, dotted, "I_{\mathfrak{A}\widehat{\otimes}\mathfrak{B}}"] \\
\mathfrak{A}\times \mathfrak{B} \arrow[r, "\widehat{\otimes}"']                                                                & \mathfrak{A}\widehat{\otimes}\mathfrak{B}.                                      
\end{tikzcd}$$ Using the 3-universal property of $\otimes$ given in lemma \ref{lem:tensorproductcategories}, the diagram can be filled by a 2-natural equivalence $\iota$, which we can also think of as an equivalence $I_{\mathfrak{A}\widehat{\otimes}\mathfrak{B}}\rightarrow I_{\mathfrak{A}}\widehat{\otimes} I_{\mathfrak{B}}$. Using this together with arguments analogous to the ones presented above, we can define the unitors, and prove that they satisfy the relevant axioms. This concludes the proof of the first part of the statement.

Now, we shall prove that the monoidal 2-functor $\widehat{\otimes}$ we have constructed above is 3-universal. Let $(\mathfrak{E}, \Box_{\mathfrak{E}}, I_{\mathfrak{E}}, \alpha_{\mathfrak{E}}, l_{\mathfrak{E}}, r_{\mathfrak{E}})$ be a monoidal $\mathds{k}$-linear 2-category that is locally Cauchy complete, and Cauchy complete. Let $$(F,\chi_F, \iota_F, \omega_F, \gamma_F,\delta_F) :\mathfrak{A}\times \mathfrak{B}\rightarrow \mathfrak{E}$$ be a monoidal bilinear 2-functor in the notation of \cite{SP}. Immediately, we see that there is a 2-natural equivalence $\phi:F'\circ \widehat{\otimes}\Rightarrow F$ for some linear 2-functor $F':\mathfrak{A}\widehat{\otimes}\mathfrak{B}\rightarrow \mathfrak{E}$. We want to show that $F'$ and $\phi$ are monoidal. To this end, consider the diagram depicted in figure \ref{fig:figure3}.

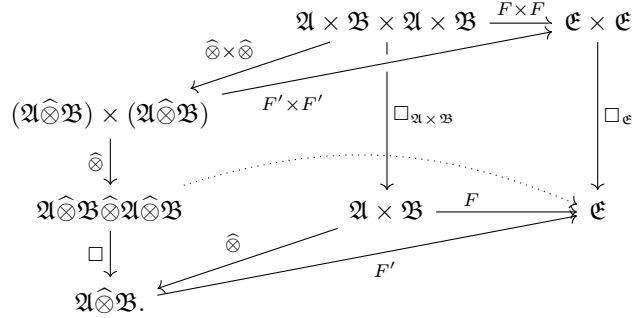
\begin{figure}[!hbt]
$$\begin{tikzcd}
                                                                                                                           & \mathfrak{A}\times\mathfrak{B}\times\mathfrak{A}\times\mathfrak{B} \arrow[r, "F\times F"] \arrow[ld, "\widehat{\otimes}\times \widehat{\otimes}"'] \arrow[dd, "\Box_{\mathfrak{A}\times\mathfrak{B}}"] & \mathfrak{E}\times\mathfrak{E} \arrow[dd, "\Box_{\mathfrak{E}}"] \\
(\mathfrak{A}\widehat{\otimes}\mathfrak{B})\times(\mathfrak{A}\widehat{\otimes}\mathfrak{B}) \arrow[rru, pos=0.1, crossing over, "F'\times F'"'] \arrow[d, "\widehat{\otimes}"'] &                                                                                                                                                                                    &                                                                  \\
\mathfrak{A}\widehat{\otimes}\mathfrak{B}\widehat{\otimes}\mathfrak{A}\widehat{\otimes}\mathfrak{B} \arrow[d, "\Box"'] \arrow[rr, dotted, bend left=20]          & \mathfrak{A}\times\mathfrak{B} \arrow[r, near start, "F"] \arrow[ld, "\widehat{\otimes}"']                                                                                                               & \mathfrak{E}                                                     \\
\mathfrak{A}\widehat{\otimes}\mathfrak{B}. \arrow[rru, "F'"']                                                                         &                                                                                                                                                                                    &                                                                 
\end{tikzcd}$$
\caption{Monoidality of $F'$}
\label{fig:figure3}
\end{figure}

The left face can be filled using $\chi$, the bottom face using $\phi$, the top face using $\phi\times\phi$, and the back face using $\chi_{F}$. Using the 3-universality of $\boxdot$ obtained in theorem \ref{thm:2Deligne}, we can factor the composite $\Box_{\mathfrak{E}}\circ (F'\times F')$ through $\boxdot$ up to a 2-natural equivalence. This fills in the upper half of the front face of the diagram represented in figure \ref{fig:figure3}, and we can then use lemma \ref{lem:multi2Deligne} to complete the bottom half, and get an invertible modification $\Pi$ witnessing the commutativity of the resulting prism. In particular, we denote the 2-natural equivalence corresponding to the composite of the front face of the diagram given in figure \ref{fig:figure3} by $\chi'_{F}$. One now checks that $\chi'_F$ endow the 2-functor $F'$ with a monoidal structure, and that $\phi$ together with $\Pi$ define a monoidal 2-natural equivalence $\phi:F\Rightarrow F'$. The remainder of the proof uses the same techniques, but involves big diagrams, which we omit.
\end{proof}

We are mainly interested in tensor 2-categories, i.e. rigid monoidal linear 2-categories (see \cite{D5}). After having proven theorem \ref{thm:CompletedTensorMon}, it is therefore natural to ask whether the completed tensor product of two tensor 2-categories is again a tensor 2-category. It turns out that the answer is positive.

\begin{Proposition}\label{prop:completedrigid}
If $\mathfrak{A}$ and $\mathfrak{B}$ be tensor 2-categories, then $\mathfrak{A}\widehat{\otimes}\mathfrak{B}$ with the monoidal structure of theorem \ref{thm:CompletedTensorMon} is a tensor 2-category.
\end{Proposition}
\begin{proof}
We prove that the objects of $\mathfrak{A}\widehat{\otimes}\mathfrak{B}$ have right duals. The existence of left duals can be treated similarly. As the 2-functor $\mathfrak{A}\times\mathfrak{B}\rightarrow\mathfrak{A}\widehat{\otimes}\mathfrak{B}$ is monoidal, every object of the form $A\widehat{\otimes}B$ for some $A$ in $\mathfrak{A}$, and $B$ in $\mathfrak{B}$ has a right dual given explicitly by $A^{\sharp}\widehat{\otimes} B^{\sharp}$. Let $E$ be any object of $\mathfrak{A}\widehat{\otimes} \mathfrak{B}$. We claim that there exists a 2-condensation monad supported on an object that is a direct sum of objects in the image of $\widehat{\otimes}$, whose splitting is $E$. Namely, by construction, we have $\widehat{\otimes} = \kappa\circ \beta$, and $\beta$ is esentially surjective. The claim now follows from the definition of the Cauchy completion as the composition of the direct sum completion and the Karoubi envelope (see appendix A.2 of \cite{D1}). Thence it is enough to prove that the splitting of a 2-condensation monad supported on an object that has a right dual has a right dual. This is precisely the content of lemma \ref{lem:dualcondensation} below.
\end{proof}

The following lemma should be compared with theorem 4.1.1 of \cite{GJF}. Because our context is different, we provide a proof in the appendix.

\begin{Lemma}\label{lem:dualcondensation}
Let $\mathfrak{C}$ be a monoidal 2-category (not necessarily linear) that is locally idempotent complete and in which every 2-condensation monad splits. Let $C$ be an object of $\mathfrak{C}$ that has a right dual, and $(C,e,\mu,\delta)$ a 2-condensation monad. Then, the splitting of $(C,e,\mu,\delta)$ has a right dual.
\end{Lemma}
\begin{proof}
See the appendix.
\end{proof}

We now assume that our base field $\mathds{k}$ is perfect. We have recalled in propositions \ref{prop:1DeligneMon} and \ref{prop:1Delignerigid} that the 1-Deligne tensor product of two finite tensor categories is itself a finite tensor category. We now prove a categorified version of this statement for finite semisimple tensor categories.

\begin{Theorem}\label{thm:2delignetensor}
Given $\mathfrak{C}$ and $\mathfrak{D}$ two compact semisimple tensor 2-categories over a perfect field. Their 2-Deligne tensor product $\mathfrak{C}\boxdot\mathfrak{D}$ is a compact semisimple tensor 2-category. Further, the 2-functor $\boxdot:\mathfrak{C}\times\mathfrak{D}\rightarrow \mathfrak{C}\boxdot\mathfrak{D}$ is monoidal.
\end{Theorem}
\begin{proof}
This follows by combining theorems \ref{thm:2Deligne} and \ref{thm:CompletedTensorMon} together with proposition \ref{prop:completedrigid}.
\end{proof}

Over an algebraically closed field or a real closed, it is shown in \cite{D5} that every compact semisimple 2-category is finite. In this case, we use the term multifusion 2-category in lieu of compact semisimple tensor 2-category. Further, a fusion 2-category is a multifusion 2-category whose monoidal unit is a simple object.

\begin{Lemma}
Let $\mathds{k}$ be an algebraically closed field. If $\mathfrak{C}$ and $\mathfrak{D}$ are two fusion 2-categories, then $\mathfrak{C}\boxdot\mathfrak{D}$ is a fusion 2-category.
\end{Lemma}
\begin{proof}
This follows from proposition \ref{prop:2DeligneHom}.
\end{proof}

We end by examining some examples.

\begin{Example}
Let us work over $\mathbb{R}$. Then, $\mathbf{Mod}(\mathbf{Vect}_{\mathbb{C}})$ is a fusion 2-category with a single equivalence class of simple object. But, we found in example \ref{ex:VectR} that $$\mathbf{Mod}(\mathbf{Vect}_{\mathbb{C}})\boxdot \mathbf{Mod}(\mathbf{Vect}_{\mathbb{C}})\simeq \mathbf{Mod}(\mathbf{Vect}_{\mathbb{C}})\boxplus \mathbf{Mod}(\mathbf{Vect}_{\mathbb{C}}).$$ In particular, the monoidal unit of this multifusion 2-category is the direct sum of the two non-equivalent simple objects of the right hand-side. Thus, $\mathbf{Mod}(\mathbf{Vect}_{\mathbb{C}})\boxdot \mathbf{Mod}(\mathbf{Vect}_{\mathbb{C}})$ is not a fusion 2-category.
\end{Example}

\begin{Example}\label{ex:2b+3triv}
We use the notation of example \ref{ex:Mod(6)}. In particular, we work over an algebraically closed field of characteristic distinct from $2$ and $3$. Let us endow $\mathcal{C}$ with a braiding via $$\mathcal{C}\simeq \mathbf{Vect}_{\mathbb{Z}/2\mathbb{Z}}^{\beta}\boxtimes \mathbf{Vect}_{\mathbb{Z}/3\mathbb{Z}}^{triv},$$ where $triv$ is the trivial symmetric braiding, and $\beta$ denotes the super braiding. The induced monoidal product on $\mathbf{Mod}(\mathcal{C})$ is given by the relative Deligne tensor product $\boxtimes_{\mathcal{C}}$. Using theorem \ref{thm:2delignetensor} and theorem \ref{thm:2Deligneessentiallysurjective}, we see that the table depicting the action of $\boxtimes_{\mathcal{C}}$ on the equivalence classes of simple objects of $\mathbf{Mod}(\mathcal{C})$ is the product of the tables for $\mathbf{Mod}(\mathbf{Vect}_{\mathbb{Z}/2\mathbb{Z}}^{\beta})$ and $\mathbf{Mod}(\mathbf{Vect}_{\mathbb{Z}/3\mathbb{Z}}^{triv})$ derived in example 3.3.7 of \cite{D5} (see also example 2.5.3 of \cite{D2}). The resulting table is given below, where we write $a:=3\in\mathbb{Z}/6\mathbb{Z}$, and $b:=2\in\mathbb{Z}/6\mathbb{Z}$. Further, we have removed all the trivial cohomology classes to make the table more legible.

\begin{center}
\begin{tabular}{ |c| c c c c| } 
 \hline
 $\boxtimes_{\mathcal{C}}$ & $\langle 0 \rangle$ & $\langle a \rangle$ & $\langle b \rangle$ & $\langle a+b \rangle$ \\ 
 \hline
 $\langle 0 \rangle$ & $\langle 0 \rangle$ & $\langle a \rangle$ & $\langle b \rangle$ & $\langle a+b \rangle$ \\ 
 $\langle a \rangle$ & $\langle a \rangle$ & $\langle 0 \rangle$ & $\langle a+b \rangle$ & $\langle b \rangle$ \\ 
 $\langle b \rangle$ & $\langle b \rangle$ & $\langle a+b \rangle$ & $3\langle b \rangle$ & $3\langle a+b \rangle$ \\ 
 $\langle a+b \rangle$ & $\langle a+b \rangle$ & $\langle b \rangle$ & $3\langle a+b \rangle$ & $3\langle b \rangle$  \\ 
 \hline
\end{tabular}
\end{center}
\end{Example}

\begin{Example}\label{ex:2b+2b}
We follow the notations of examples \ref{ex:Mod(2+2)}. In particular, we work over an algebraically closed field of characteristic not equal to $2$. The equivalence of fusion categories $$\mathcal{D}\simeq \mathbf{Vect}_{\mathbb{Z}/2\mathbb{Z}}^{\beta}\boxtimes \mathbf{Vect}_{\mathbb{Z}/2\mathbb{Z}}^{\beta}$$ induces a braiding on $\mathcal{D}$. The induced monoidal structure on $\mathbf{Mod}(\mathcal{D})$ is given by $\boxtimes_{\mathcal{D}}$. Using theorem \ref{thm:2delignetensor} we can partially compute the action of $\boxtimes_{\mathcal{D}}$ on the equivalence classes of simple objects of $\mathbf{Mod}(\mathcal{D})$. Namely, we can compute the monoidal products (depicted in the table below) involving the objects $\langle 0\rangle$, $\langle a\rangle$, $\langle b\rangle$, and $\langle a,b\rangle$. We have omitted the trivial cohomology classes to declutter the table, and denoted the simple object $(\langle a,b \rangle, \nu)$ succinctly by $\nu$.

\begin{center}
\begin{tabular}{ |c| c c c c c c| } 
 \hline
 $\boxtimes_{\mathcal{D}}$ & $\langle 0 \rangle$ & $\langle a \rangle$ & $\langle b \rangle$ & $\langle c \rangle$ & $\langle a,b \rangle$ & $\nu$ \\ 
 \hline
 $\langle 0 \rangle$ & $\langle 0 \rangle$ & $\langle a \rangle$ & $\langle b \rangle$ & $\langle c \rangle$ & $\langle a,b \rangle$ & $\nu$ \\ 
 $\langle a \rangle$ & $\langle a \rangle$ & $\langle 0 \rangle$ & $\langle a,b \rangle$ & $\nu$ & $\langle b \rangle$ & $\langle c\rangle$ \\ 
 $\langle b \rangle$ & $\langle b \rangle$ & $\langle a,b \rangle$ & $\langle 0 \rangle$ & $\nu$ &  $\langle a \rangle$ & $\langle c\rangle$ \\ 
 $\langle c \rangle$ & $\langle c \rangle$ & $\nu$ & $\nu$ & $2\langle c \rangle$ & $\langle c\rangle$ & $2\nu$ \\ 
 $\langle a,b \rangle$ & $\langle a,b \rangle$ & $\langle b \rangle$ & $\langle a \rangle$ & $\langle c\rangle$ & $\langle 0 \rangle$ & $\nu$  \\ 
 $\nu$ & $\nu$ & $\langle c\rangle$ & $\langle c\rangle$ & $2\nu$ & $\nu$ & $2\langle c\rangle$  \\ 
 \hline
\end{tabular}
\end{center}

\noindent Provided our base field has characteristic zero, the remaining entries can be filled using proposition 2.5.1 of \cite{D2}. Over an algebraically closed field of characteristic greater than 2, they can be computed using a similar method.
\end{Example}

\begin{Remark}
As the fusion categories $\mathcal{C}$ and $\mathcal{D}$ of examples \ref{ex:2b+3triv} and \ref{ex:2b+2b} are symmetric, so are $\mathbf{Mod}(\mathcal{C})$ and $\mathbf{Mod}(\mathcal{D})$. This explains why their respective multiplication tables are symmetric along the diagonal.
\end{Remark}

\appendix

\section*{Appendix}
\addcontentsline{toc}{section}{Appendix}

\renewcommand*{\proofname}{Proof of lemma \ref{lem:dualcondensation}}

\begin{proof}
It is convenient to assume that $\mathfrak{C}$ is strict cubical.Following the notation of \cite{SP}, we write $\phi$ for the interchanger. Observe that, even though $\mathfrak{C}$ is not a priori rigid, it makes sense to consider the right dual 2-functor $(-)^{\sharp}$ constructed in lemma A.2.2 of \cite{D2} on the full sub-2-category on the objects of $\mathfrak{C}$ that have a right dual.

The image of the 2-condensation monad $M:=(C,e,\mu,\delta)$ in $\mathfrak{C}$ under $(-)^{\sharp}$ yields a 2-condensation monad $N:=(C^{\sharp},e^{\sharp}, \mu^{\sharp}, \delta^{\sharp})$ in $\mathfrak{C}^{1op}$, which we will treat as a 2-condensation monad in $\mathfrak{C}$. We write $(C,D,f,g,\gamma,\psi)$ and $\theta:g\circ f\cong e$ for a splitting of $M$ in $\mathfrak{C}$, as well as $(C^{\sharp},D',f^{\sharp},g',\gamma',\psi')$ and $\theta':g'\circ f^{\sharp}\cong e^{\sharp}$ for a splitting of $N$. Further, recall that splittings of 2-condensation moands are preserved by all 2-functors. In particular, the splitting for the monoidal product of two 2-condensation monads is given by the monoidal product of their splittings.

We claim that $D'$ is a right dual for $D$. Because $C$ has a right dual $C^{\sharp}$, there exists 1-morphisms $i_C:I\rightarrow C^{\sharp}\Box C$ and $C\Box C^{\sharp}\rightarrow I$ satisfying the snake equations up to 2-isomorphisms. We fix such a 2-isomorphism $$E_C:(C^{\sharp}\Box e_C)\circ (i_C\Box C^{\sharp})\cong Id_{C^{\sharp}}\textrm{ and }F_C:Id_C\cong (e_C\Box C)\circ (C\Box i_C).$$

Observe that the 1-morphisms in $\mathfrak{C}$ $$i_M:= (C^{\sharp}\Box e)\circ i_C\ \ \text{ and }\ \ e_M:= e_C \circ (e\Box C^{\sharp})$$ can canonically be upgraded to bimodules (in the sense of \cite{D1} definition A.2.1, see also \cite{GJF}). More precisely, if we abuse the notation, and use $I$ to denote the trivial 2-condensation monad on $I$, $i_M$ is an $(N\Box M, I)$-bimodule, and $e_M$ is an $(I,M\Box N)$-bimodule. For the sake brevity, we shall only explicitly write down the two structure 2-morphisms that will be relevant for our purposes, namely $\nu^l_{i_M}$ and $\beta^r_{e_M}$. These 2-morphisms are depicted below using the wire diagram calculus of \cite{B}:

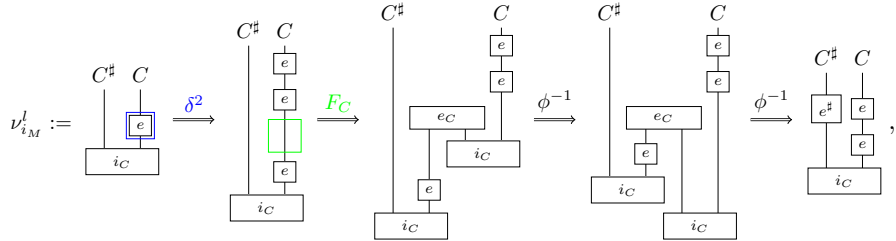
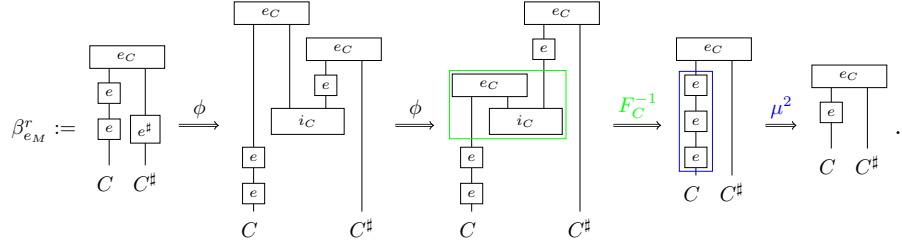
$$\adjustbox{scale =0.8}{$\nu^l_{i_M}:=
\begin{tz}[xscale=3, yscale=1.5]
\node (1) at (0,0)
{
$\begin{tz}[xscale=0.6, yscale=0.6]
  \draw (0,0.5)  to (0,2.5) node[above] {$C^{\sharp}$};
  \draw (1,0.5) to (1, 1.5) node[morphlabel] (r) {$e$} to (1,2.5)node[above] {$C$};
  
\node[morphlabel] at (0.5, 0.5) {$\,\,\,\, \,\,\, i_C \,\,\,\, \,\,\,$};

\draw[blue] ([xshift=-2pt, yshift=2pt] r.north west) rectangle ([xshift=2pt, yshift=-2pt] r.south east);
\end{tz}$
};

\node (2) at (0.8,0)
{
$\begin{tz}[xscale=0.6, yscale=0.6]
  \draw (0,0.5)  to (0,5) node[above] {$C^{\sharp}$};
  \draw (1,0.5) to (1,1.5) node[morphlabel] {$e$} to (1, 2.5) node (r) {} to (1,3.5) node[morphlabel] {$e$} to (1,4.5) node[morphlabel] {$e$} to (1,5) node[above] {$C$};
  
\node[morphlabel] at (0.5, 0.5) {$\,\,\,\, \,\,\, i_C \,\,\,\, \,\,\,$};

\draw[green] ([xshift=-7pt, yshift=7pt] r.north west) rectangle ([xshift=7pt, yshift=-7pt] r.south east);
\end{tz}$
};

\node (3) at (1.8,0)
{
$\begin{tz}[xscale=0.6, yscale=0.6]
  \draw (0,0.5)  to (0,6) node[above] {$C^{\sharp}$};
  \draw (1,0.5) to (1,1.5) node[morphlabel] {$e$} to (1, 3.5);
  \draw (2,2.5) to (2,3.5);
  \draw (3,2.5) to (3,4.5) node[morphlabel] {$e$} to (3,5.5) node[morphlabel] {$e$} to (3,6) node[above] {$C$};
  
\node[morphlabel] at (0.5, 0.5) {$\,\,\,\, \,\,\, i_C \,\,\,\, \,\,\,$};

\node[morphlabel] at (2.5, 2.5) {$\,\,\,\, \,\,\, i_C \,\,\,\, \,\,\,$};

\node[morphlabel] at (1.5, 3.5) {$\,\,\,\, \,\,\, e_C \,\,\,\, \,\,\,$};

\end{tz}$
};

\node (4) at (3,0)
{
$\begin{tz}[xscale=0.6, yscale=0.6]
  \draw (0,1.5)  to (0,6) node[above] {$C^{\sharp}$};
  \draw (1,1.5) to (1,2.5) node[morphlabel] {$e$} to (1, 3.5);
  \draw (2,0.5) to (2,3.5);
  \draw (3,0.5) to (3,4.5) node[morphlabel] {$e$} to (3,5.5) node[morphlabel] {$e$} to (3,6) node[above] {$C$};
  
\node[morphlabel] at (0.5, 1.5) {$\,\,\,\, \,\,\, i_C \,\,\,\, \,\,\,$};

\node[morphlabel] at (2.5, 0.5) {$\,\,\,\, \,\,\, i_C \,\,\,\, \,\,\,$};

\node[morphlabel] at (1.5, 3.5) {$\,\,\,\, \,\,\, e_C \,\,\,\, \,\,\,$};

\end{tz}$
};

\node (5) at (4,0)
{
$\begin{tz}[xscale=0.6, yscale=0.6]
  \draw (0,0.5) to (0,2.5) node[morphlabel] {$e^{\sharp}$}  to (0,3.5) node[above] {$C^{\sharp}$};
  \draw (1,0.5) to (1,1.5) node[morphlabel] {$e$} to (1, 2.5) node[morphlabel] {$e$} to (1,3.5) node[above] {$C$};
  
\node[morphlabel] at (0.5, 0.5) {$\,\,\,\, \,\,\, i_C \,\,\,\, \,\,\,$};
\end{tz}$
};
 
\begin{scope}[double arrow scope]
	\draw ([xshift=0.1em] 1.east) -- node[above, blue] {$\delta^2$} ([xshift=-0.1em] 2.west);
\end{scope}

\begin{scope}[double arrow scope]
	\draw ([xshift=0.1em] 2.east) -- node[above, green] {$F_C$} ([xshift=-0.1em] 3.west);
\end{scope}

\begin{scope}[double arrow scope]
	\draw ([xshift=0.1em] 3.east) -- node[above] {$\phi^{-1}$} ([xshift=-0.1em] 4.west);
\end{scope}

\begin{scope}[double arrow scope]
	\draw ([xshift=0.1em] 4.east) -- node[above] {$\phi^{-1}$} ([xshift=-0.1em] 5.west);
\end{scope}
\end{tz}$},$$

$$\adjustbox{scale =0.8}{
$\beta^r_{e_M}:=
\begin{tz}[xscale=3, yscale=1.5]
\node (1) at (0,0)
{
$\begin{tz}[xscale=0.6, yscale=0.6]
  \draw (0,0.5) node[below] {$C\phantom{^{\sharp}}$} to (0,1.5) node[morphlabel] {$e$} to (0, 2.5) node[morphlabel] {$e$}  to (0,3.5);
  \draw (1,0.5) node[below] {$C^{\sharp}$} to (1,1.5) node[morphlabel] {$e^{\sharp}$} to (1,3.5);
  
\node[morphlabel] at (0.5, 3.5) {$\,\,\,\, \,\,\, e_C \,\,\,\, \,\,\,$};
\end{tz}$
};

\node (2) at (1,0)
{
$\begin{tz}[xscale=0.6, yscale=0.6]
  \draw (0,0) node[below] {$C\phantom{^{\sharp}}$} to (0,0.5) node[morphlabel] {$e$} to (0, 1.5) node[morphlabel] {$e$} to (0,5.5);
  \draw (1, 2.5) to (1,5.5);
  \draw (2,2.5) to (2,3.5) node[morphlabel] {$e$} to (2,4.5);
  \draw (3,0) node[below] {$C^{\sharp}$} to (3,4.5);
  
\node[morphlabel] at (0.5, 5.5) {$\,\,\,\, \,\,\, e_C \,\,\,\, \,\,\,$};

\node[morphlabel] at (2.5, 4.5) {$\,\,\,\, \,\,\, e_C \,\,\,\, \,\,\,$};

\node[morphlabel] at (1.5, 2.5) {$\,\,\,\, \,\,\, i_C \,\,\,\, \,\,\,$};
\end{tz}$
};

\node (3) at (2.2,0)
{
$\begin{tz}[xscale=0.6, yscale=0.6]
  \draw (0,0) node[below] {$C\phantom{^{\sharp}}$} to (0,0.5) node[morphlabel] {$e$} to (0, 1.5) node[morphlabel] {$e$} to (0,3.5);
  \draw (1, 2.5) to (1,3) node (m) {} to (1,3.5);
  \draw (2,2.5) to (2,4.5) node[morphlabel] {$e$} to (2,5.5);
  \draw (3,0) node[below] {$C^{\sharp}$} to (3,5.5);
  
\node[morphlabel] at (0.5, 3.5) {$\,\,\,\, \,\,\, e_C \,\,\,\, \,\,\,$};
\node[morphlabel] at (2.5, 5.5) {$\,\,\,\, \,\,\, e_C \,\,\,\, \,\,\,$};
\node[morphlabel] at (1.5, 2.5) {$\,\,\,\, \,\,\, i_C \,\,\,\, \,\,\,$};

\draw[green] ([xshift=-40pt, yshift=20pt] m.north west) rectangle ([xshift=40pt, yshift=-22pt] m.south east);
\end{tz}$
};

\node (4) at (3.25,0)
{
$\begin{tz}[xscale=0.6, yscale=0.6]
  \draw (0,0) node[below] {$C\phantom{^{\sharp}}$} to (0, 0.5) node[morphlabel] {$e$} to (0, 1.5) node[morphlabel] (l) {$e$} to (0, 2.5) node[morphlabel] {$e$} to (0,3.5) ;
  \draw (1,0) node[below] {$C^{\sharp}$} to (1,3.5);
  
\node[morphlabel] at (0.5, 3.5) {$\,\,\,\, \,\,\, e_C \,\,\,\, \,\,\,$};

\draw[blue] ([xshift=-4pt, yshift=31pt] l.north west) rectangle ([xshift=4pt, yshift=-31pt] l.south east);
\end{tz}$
};

\node (5) at (4,0)
{
$\begin{tz}[xscale=0.6, yscale=0.6]
  \draw (0,0.5) node[below] {$C\phantom{^{\sharp}}$} to (0, 1.5) node[morphlabel] (r) {$e$} to (0,2.5) ;
  \draw (1,0.5) node[below] {$C^{\sharp}$} to (1,2.5);
  
\node[morphlabel] at (0.5, 2.5) {$\,\,\,\, \,\,\, e_C \,\,\,\, \,\,\,$};

\end{tz}$
};
 
\begin{scope}[double arrow scope]
	\draw ([xshift=0.1em] 1.east) -- node[above] {$\phi$} ([xshift=-0.1em] 2.west);
\end{scope}

\begin{scope}[double arrow scope]
	\draw ([xshift=0.1em] 2.east) -- node[above] {$\phi$} ([xshift=-0.1em] 3.west);
\end{scope}

\begin{scope}[double arrow scope]
	\draw ([xshift=0.1em] 3.east) -- node[above,green] {$F_C^{-1}$} ([xshift=-0.1em] 4.west);
\end{scope}

\begin{scope}[double arrow scope]
	\draw ([xshift=0.1em] 4.east) -- node[above,blue] {$\mu^2$} ([xshift=-0.1em] 5.west);
\end{scope}

\end{tz}$}.$$

Recall the construction introduced in \cite{GJF} and extensively used in appendix A.2 of \cite{D1} by which a bimodule between split 2-condensation monads yields a 1-morphism between their splittings. For instance, if we take as input the bimodule $i_M$, we use the fact that the composite 1-morphism $(f^{\sharp}\Box f)\circ i_M: I\rightarrow C^{\sharp}\Box C$ in $\mathfrak{C}$ supports a canonical idempotent given by $$\adjustbox{scale =0.8}{$
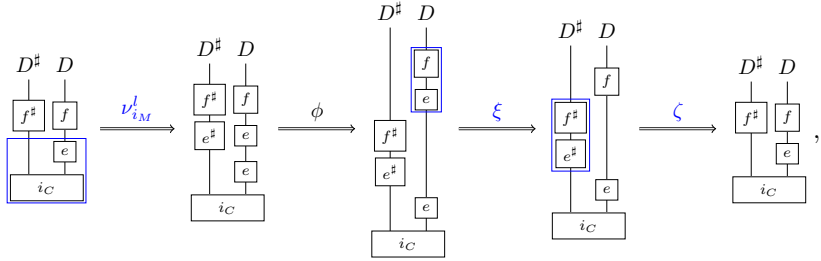
\begin{tz}[xscale=3, yscale=1.5]
\node (1) at (0,0)
{
$\begin{tz}[xscale=0.6, yscale=0.6]
  \draw (0,0.5) to (0,2.5) node[morphlabel]{$f^{\sharp}$} to (0,3.5) node[above] {$D^{\sharp}$};
  \draw (1,0.5) to (1, 1.5) node[morphlabel]{$e$} to (1, 2.5) node[morphlabel] {$f$} to (1,3.5)node[above] {$D$};
  
\node[morphlabel] at (0.5, 0.5) {$\,\,\,\, \,\,\, i_C \,\,\,\, \,\,\,$};

\draw[blue] ([xshift=-8pt, yshift=2pt] r.north west) rectangle ([xshift=36pt, yshift=-32pt] r.south east);
\end{tz}$
};

\node (2) at (1,0)
{
$\begin{tz}[xscale=0.6, yscale=0.6]
  \draw (0,0.5) to (0,2.5) node[morphlabel] {$e^{\sharp}$} to (0,3.5) node[morphlabel] {$f^{\sharp}$} to (0,4.5) node[above] {$D^{\sharp}$};
  \draw (1,0.5) to (1,1.5) node[morphlabel] {$e$} to (1, 2.5) node[morphlabel] {$e$} to (1,3.5) node[morphlabel] {$f$} to (1,4.5) node[above] {$D$};
  
\node[morphlabel] at (0.5, 0.5) {$\,\,\,\, \,\,\, i_C \,\,\,\, \,\,\,$};
\end{tz}$
};

\node (3) at (2,0)
{
$\begin{tz}[xscale=0.6, yscale=0.6]
  \draw (0,0.5) to (0,2.5) node[morphlabel] {$e^{\sharp}$} to (0,3.5) node[morphlabel] {$f^{\sharp}$} to (0,6.5) node[above] {$D^{\sharp}$};
  \draw (1,0.5) to (1,1.5) node[morphlabel] {$e$} to (1, 4.5) node[morphlabel] {$e$} to (1,5) node(r) {} to (1,5.5) node[morphlabel] {$f$} to (1,6.5) node[above] {$D$};
  
\node[morphlabel] at (0.5, 0.5) {$\,\,\,\, \,\,\, i_C \,\,\,\, \,\,\,$};

\draw[blue] ([xshift=-6pt, yshift=21pt] r.north west) rectangle ([xshift=6pt, yshift=-18pt] r.south east);
\end{tz}$
};

\node (4) at (3,0)
{
$\begin{tz}[xscale=0.6, yscale=0.6]
  \draw (0,0.5) to (0,2.5) node[morphlabel] {$e^{\sharp}$} to (0,3) node(r) {} to (0,3.5) node[morphlabel] {$f^{\sharp}$} to (0,5.5) node[above] {$D^{\sharp}$};
  \draw (1,0.5) to (1,1.5) node[morphlabel] {$e$} to (1, 4.5) node[morphlabel] {$f$} to (1,5.5) node[above] {$D$};
  
\node[morphlabel] at (0.5, 0.5) {$\,\,\,\, \,\,\, i_C \,\,\,\, \,\,\,$};
\draw[blue] ([xshift=-9pt, yshift=23pt] r.north west) rectangle ([xshift=9pt, yshift=-22pt] r.south east);
\end{tz}$
};

\node (5) at (4,0)
{
$\begin{tz}[xscale=0.6, yscale=0.6]
  \draw (0,0.5) to (0,2.5) node[morphlabel]{$f^{\sharp}$} to (0,3.5) node[above] {$D^{\sharp}$};
  \draw (1,0.5) to (1, 1.5) node[morphlabel]{$e$} to (1, 2.5) node[morphlabel] {$f$} to (1,3.5)node[above] {$D$};
  
\node[morphlabel] at (0.5, 0.5) {$\,\,\,\, \,\,\, i_C \,\,\,\, \,\,\,$};

\end{tz}$
};
 
\begin{scope}[double arrow scope]
	\draw ([xshift=0.1em] 1.east) -- node[above, blue] {$\nu^l_{i_M}$} ([xshift=-0.1em] 2.west);
\end{scope}

\begin{scope}[double arrow scope]
	\draw ([xshift=0.1em] 2.east) -- node[above] {$\phi$} ([xshift=-0.1em] 3.west);
\end{scope}

\begin{scope}[double arrow scope]
	\draw ([xshift=0.1em] 3.east) -- node[above, blue] {$\xi$} ([xshift=-0.1em] 4.west);
\end{scope}

\begin{scope}[double arrow scope]
	\draw ([xshift=0.1em] 4.east) -- node[above, blue] {$\zeta$} ([xshift=-0.1em] 5.west);
\end{scope}
\end{tz}$},$$ where $\xi=(\phi\circ f)\cdot (f\circ \theta^{-1})$ and $\zeta=(\phi'\circ f^{\sharp})\cdot (f^{\sharp}\circ \theta'^{-1})$. The splitting of this canonical idempotent gives us the desired 1-morphism $i_D:I\rightarrow D'\Box D$. Proceeding analogously with the bimodule $e_M$, we get a 1-morphism $e_D:D'\Box D\rightarrow I$. We claim that these two 1-morphisms satisfy the snake equations up to 2-isomorphism. To this end, consider the composite $(M,M)$-bimodule $$(e_M\Box e)\otimes_{M\Box N\Box M} (e\Box i_M).$$ Because the construction reviewed above sends tensor products of bimodules to composites of 1-morphisms (see proposition 1.2.3 of \cite{D1}), we have that $(e_D\Box D)\circ (D\Box i_D)$ is a splitting of the canonical idempotent supported on $f\circ (e_M\Box e)\otimes_{M\Box N\Box M} (e\Box i_M)\circ g$. Thus, it only remains to prove that $Id_D$ is a splitting for this idempotent. This is readily implied by the claim that as $(M,M)$-bimodules, we have $$e \cong (e_M\Box e)\otimes_{M\Box N\Box M} (e\Box i_M).$$ Namely, the splitting of the canonical idempotent supported on the $(M,M)$-bimodule $e$ is $Id_D$. Now, the $(M,M)$-bimodule $(e_M\Box e)\otimes_{M\Box N\Box M} (e\Box i_M)$ is defined as the splitting of the idempotent $\tau$ supported on $(e_M\Box e)\circ (e\Box i_M)$ given by $$\adjustbox{scale =0.8}{$
\tau:=
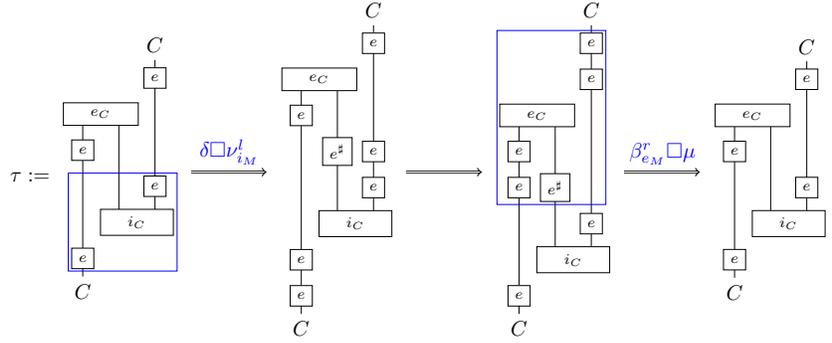
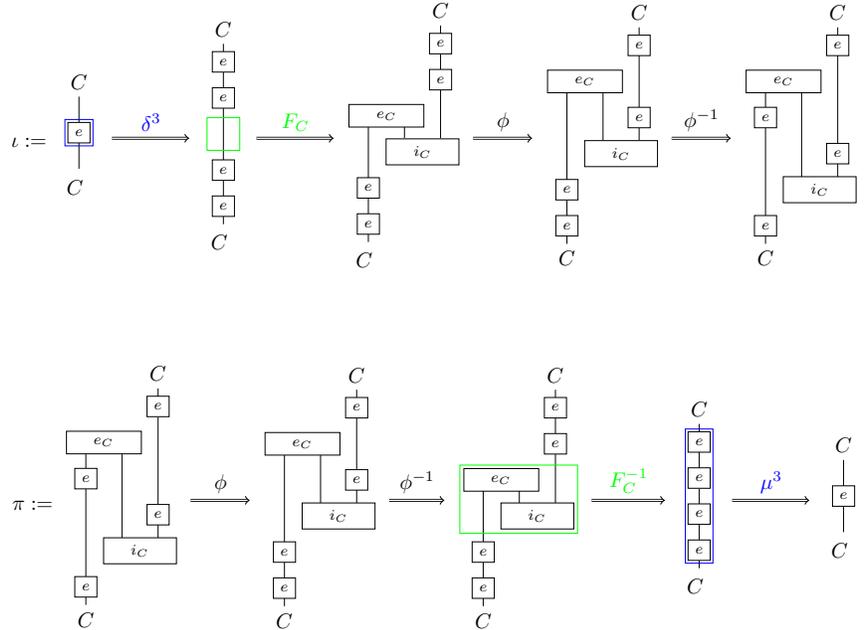
\begin{tz}[xscale=3, yscale=1.5]
\node (1) at (0,0)
{
$\begin{tz}[xscale=0.6, yscale=0.6]
  \draw (0,0) node[below] {$C$} to (0,0.5) node[morphlabel] {$e$} to (0, 3.5) node[morphlabel] {$e$} to (0,4.5);
  \draw (1, 1.5) to (1,4.5);
  \draw (2,1.5) to (2,2.5) node[morphlabel] {$e$} to (2,5.5) node[morphlabel] {$e$} to (2,6) node[above] {$C$};
  
\node[morphlabel] at (0.5, 4.5) {$\,\,\,\, \,\,\, e_C \,\,\,\, \,\,\,$};

\node[morphlabel] at (1.5, 1.5) {$\,\,\,\, \,\,\, i_C \,\,\,\, \,\,\,$};

\draw[blue] ([xshift=-34 pt, yshift=-87pt] m.north west) rectangle ([xshift=40pt, yshift=2pt] m.south east);
\end{tz}$
};

\node (2) at (1.2,0)
{
$\begin{tz}[xscale=0.6, yscale=0.6]
  \draw (0,0) node[below] {$C$} to (0,0.5) node[morphlabel] {$e$} to (0,1.5) node[morphlabel] {$e$} to (0, 5.5) node[morphlabel] {$e$} to (0,6.5);
  \draw (1, 2.5) to (1, 4.5) node[morphlabel] {$e^{\sharp}$} to (1,6.5);
  \draw (2,2.5) to (2,3.5) node[morphlabel] {$e$} to (2, 4.5) node[morphlabel] {$e$} to (2,7.5) node[morphlabel] {$e$} to (2,8) node[above] {$C$};
  
\node[morphlabel] at (0.5, 6.5) {$\,\,\,\, \,\,\, e_C \,\,\,\, \,\,\,$};

\node[morphlabel] at (1.5, 2.5) {$\,\,\,\, \,\,\, i_C \,\,\,\, \,\,\,$};
\end{tz}$
};

\node (3) at (2.4,0){
$\begin{tz}[xscale=0.6, yscale=0.6]
  \draw (0,0) node[below] {$C$} to (0,0.5) node[morphlabel] {$e$} to (0,3.5) node[morphlabel] {$e$} to (0, 4.5) node[morphlabel] {$e$} to (0,5.5);
  \draw (1, 1.5) to (1, 3.5) node[morphlabel] {$e^{\sharp}$} to (1,5.5);
  \draw (2,1.5) to (2,2.5) node[morphlabel] {$e$} to (2, 6.5) node[morphlabel] {$e$} to (2,7.5) node[morphlabel] {$e$} to (2,8) node[above] {$C$};
  
\node[morphlabel] at (0.5, 5.5) {$\,\,\,\, \,\,\, e_C \,\,\,\, \,\,\,$};

\node[morphlabel] at (1.5, 1.5) {$\,\,\,\, \,\,\, i_C \,\,\,\, \,\,\,$};

\draw[blue] ([xshift=-40 pt, yshift=-5pt] m.north west) rectangle ([xshift=34pt, yshift=144pt] m.south east);
\end{tz}$
};

\node (4) at (3.6,0)
{
$\begin{tz}[xscale=0.6, yscale=0.6]
  \draw (0,0) node[below] {$C$} to (0,0.5) node[morphlabel] {$e$} to (0, 3.5) node[morphlabel] {$e$} to (0,4.5);
  \draw (1, 1.5) to (1,4.5);
  \draw (2,1.5) to (2,2.5) node[morphlabel] {$e$} to (2,5.5) node[morphlabel] {$e$} to (2,6) node[above] {$C$};
  
\node[morphlabel] at (0.5, 4.5) {$\,\,\,\, \,\,\, e_C \,\,\,\, \,\,\,$};

\node[morphlabel] at (1.5, 1.5) {$\,\,\,\, \,\,\, i_C \,\,\,\, \,\,\,$};
\end{tz}$
};
 
\begin{scope}[double arrow scope]
	\draw ([xshift=0.1em] 1.east) -- node[above, blue] {$\delta\Box \nu^l_{i_M}$} ([xshift=-0.1em] 2.west);
\end{scope}

\begin{scope}[double arrow scope]
	\draw ([xshift=0.1em] 2.east) -- node[above] {} ([xshift=-0.1em] 3.west);
\end{scope}

\begin{scope}[double arrow scope]
	\draw ([xshift=0.1em] 3.east) -- node[above, blue] {$\beta^r_{e_M}\Box \mu$} ([xshift=-0.1em] 4.west);
\end{scope}
\end{tz}
$}.$$ The 2-morphisms

$$\adjustbox{scale =0.8}{$
\iota:=
\begin{tz}[xscale=3, yscale=1.5]
\node (1) at (0,0)
{
$\begin{tz}[xscale=0.6, yscale=0.6]
  \draw (0,0.5) node[below] {$C\phantom{^{\sharp}}$} to (0, 1.5) node[morphlabel] (m) {$e$} to (0,2.5) node[above] {$C$};
  
  \draw[blue] ([xshift=-2pt, yshift=2pt] m.north west) rectangle ([xshift=2pt, yshift=-2pt] m.south east);
\end{tz}$
};

\node (2) at (0.8,0)
{
$\begin{tz}[xscale=0.6, yscale=0.6]
  \draw (0,0) node[below] {$C\phantom{^{\sharp}}$} to (0, 0.5) node[morphlabel] {$e$} to (0, 1.5) node[morphlabel] {$e$} to (0,2.5) node (m) {} to (0, 3.5) node[morphlabel] {$e$} to (0,4.5) node[morphlabel] {$e$} to (0,5) node[above] {$C$};

\draw[green] ([xshift=-7pt, yshift=7pt] m.north west) rectangle ([xshift=7pt, yshift=-7pt] m.south east);
\end{tz}$
};

\node (3) at (1.8,0)
{
$\begin{tz}[xscale=0.6, yscale=0.6]
  \draw (0,0) node[below] {$C\phantom{^{\sharp}}$} to (0, 0.5) node[morphlabel] {$e$} to (0, 1.5) node[morphlabel] {$e$} to (0,3.5);
  \draw (1, 2.5) to (1,3.5);
  \draw (2,2.5) to (2, 4.5) node[morphlabel] {$e$} to (2,5.5) node[morphlabel] {$e$} to (2,6) node[above] {$C$};
  
  \node[morphlabel] at (0.5, 3.5) {$\,\,\,\, \,\,\, e_C \,\,\,\, \,\,\,$};
  
  \node[morphlabel] at (1.5, 2.5) {$\,\,\,\, \,\,\, i_C \,\,\,\, \,\,\,$};
\end{tz}$
};
 
\node (4) at (2.9,0)
{
$\begin{tz}[xscale=0.6, yscale=0.6]
  \draw (0,0) node[below] {$C$} to (0,0.5) node[morphlabel] {$e$} to (0, 1.5) node[morphlabel] {$e$} to (0,4.5);
  \draw (1, 2.5) to (1,4.5);
  \draw (2,2.5) to (2,3.5) node[morphlabel] {$e$} to (2,5.5) node[morphlabel] {$e$} to (2,6) node[above] {$C$};
  
\node[morphlabel] at (0.5, 4.5) {$\,\,\,\, \,\,\, e_C \,\,\,\, \,\,\,$};

\node[morphlabel] at (1.5, 2.5) {$\,\,\,\, \,\,\, i_C \,\,\,\, \,\,\,$};
\end{tz}$
};

\node (5) at (4,0)
{
$\begin{tz}[xscale=0.6, yscale=0.6]
  \draw (0,0) node[below] {$C$} to (0,0.5) node[morphlabel] {$e$} to (0, 3.5) node[morphlabel] {$e$} to (0,4.5);
  \draw (1, 1.5) to (1,4.5);
  \draw (2,1.5) to (2,2.5) node[morphlabel] {$e$} to (2,5.5) node[morphlabel] {$e$} to (2,6) node[above] {$C$};
  
\node[morphlabel] at (0.5, 4.5) {$\,\,\,\, \,\,\, e_C \,\,\,\, \,\,\,$};

\node[morphlabel] at (1.5, 1.5) {$\,\,\,\, \,\,\, i_C \,\,\,\, \,\,\,$};
\end{tz}$
};
 
\begin{scope}[double arrow scope]
	\draw ([xshift=0.1em] 1.east) -- node[above, blue] {$\delta^3$} ([xshift=-0.1em] 2.west);
\end{scope}

\begin{scope}[double arrow scope]
	\draw ([xshift=0.1em] 2.east) -- node[above,green] {$F_C$} ([xshift=-0.1em] 3.west);
\end{scope}

\begin{scope}[double arrow scope]
	\draw ([xshift=0.1em] 3.east) -- node[above] {$\phi$} ([xshift=-0.1em] 4.west);
\end{scope}

\begin{scope}[double arrow scope]
	\draw ([xshift=0.1em] 4.east) -- node[above] {$\phi^{-1}$} ([xshift=-0.1em] 5.west);
\end{scope}
\end{tz}
$},$$

$$\adjustbox{scale =0.8}{$
\pi:=
\begin{tz}[xscale=3, yscale=1.5]
\node (1) at (0,0)
{
$\begin{tz}[xscale=0.6, yscale=0.6]
  \draw (0,0) node[below] {$C$} to (0,0.5) node[morphlabel] {$e$} to (0, 3.5) node[morphlabel] {$e$} to (0,4.5);
  \draw (1, 1.5) to (1,4.5);
  \draw (2,1.5) to (2,2.5) node[morphlabel] {$e$} to (2,5.5) node[morphlabel] {$e$} to (2,6) node[above] {$C$};
  
\node[morphlabel] at (0.5, 4.5) {$\,\,\,\, \,\,\, e_C \,\,\,\, \,\,\,$};

\node[morphlabel] at (1.5, 1.5) {$\,\,\,\, \,\,\, i_C \,\,\,\, \,\,\,$};
\end{tz}$
};

\node (2) at (1.1,0)
{
$\begin{tz}[xscale=0.6, yscale=0.6]
  \draw (0,0) node[below] {$C$} to (0,0.5) node[morphlabel] {$e$} to (0, 1.5) node[morphlabel] {$e$} to (0,4.5);
  \draw (1, 2.5) to (1,4.5);
  \draw (2,2.5) to (2,3.5) node[morphlabel] {$e$} to (2,5.5) node[morphlabel] {$e$} to (2,6) node[above] {$C$};
  
\node[morphlabel] at (0.5, 4.5) {$\,\,\,\, \,\,\, e_C \,\,\,\, \,\,\,$};

\node[morphlabel] at (1.5, 2.5) {$\,\,\,\, \,\,\, i_C \,\,\,\, \,\,\,$};
\end{tz}$
};

\node (3) at (2.2,0){
$\begin{tz}[xscale=0.6, yscale=0.6]
  \draw (0,0) node[below] {$C$} to (0,0.5) node[morphlabel] {$e$} to (0, 1.5) node[morphlabel] {$e$} to (0,3.5);
  \draw (1, 2.5) to (1,3) node (m) {} to (1,3.5);
  \draw (2,2.5) to (2,4.5) node[morphlabel] {$e$} to (2,5.5) node[morphlabel] {$e$} to (2,6) node[above] {$C$};
  
\node[morphlabel] at (0.5, 3.5) {$\,\,\,\, \,\,\, e_C \,\,\,\, \,\,\,$};

\node[morphlabel] at (1.5, 2.5) {$\,\,\,\, \,\,\, i_C \,\,\,\, \,\,\,$};

\draw[green] ([xshift=-41pt, yshift=20pt] m.north west) rectangle ([xshift=40pt, yshift=-22pt] m.south east);
\end{tz}$
};

\node (4) at (3.2,0)
{
$\begin{tz}[xscale=0.6, yscale=0.6]
  \draw (0,0) node[below] {$C\phantom{^{\sharp}}$} to (0, 0.5) node[morphlabel] {$e$} to (0, 1.5) node[morphlabel] {$e$} to (0,2) node (m) {} to (0, 2.5) node[morphlabel] {$e$} to (0,3.5) node[morphlabel] {$e$} to (0,4) node[above] {$C$};

\draw[blue] ([xshift=-5pt, yshift=47pt] m.north west) rectangle ([xshift=5pt, yshift=-47pt] m.south east);
\end{tz}$
};

\node (5) at (4,0)
{
$\begin{tz}[xscale=0.6, yscale=0.6]
  \draw (0,0.5) node[below] {$C\phantom{^{\sharp}}$} to (0, 1.5) node[morphlabel] {$e$} to (0,2.5) node[above] {$C$};
\end{tz}$
};
 
\begin{scope}[double arrow scope]
	\draw ([xshift=0.1em] 1.east) -- node[above] {$\phi$} ([xshift=-0.1em] 2.west);
\end{scope}

\begin{scope}[double arrow scope]
	\draw ([xshift=0.1em] 2.east) -- node[above] {$\phi^{-1}$} ([xshift=-0.1em] 3.west);
\end{scope}

\begin{scope}[double arrow scope]
	\draw ([xshift=0.1em] 3.east) -- node[above,green] {$F_C^{-1}$} ([xshift=-0.1em] 4.west);
\end{scope}

\begin{scope}[double arrow scope]
	\draw ([xshift=0.1em] 4.east) -- node[above,blue] {$\mu^3$} ([xshift=-0.1em] 5.west);
\end{scope}
\end{tz}
$}$$

\noindent provide a splitting for the idempotent $\tau$ in the category of $(M,M)$-bimodules, which proves the claim. Checking that the remaining snake equation holds up to 2-isomorphism can be done using analogous techniques. We leave the details to the reader.
\end{proof}

\bibliography{bibliography.bib}

\end{document}